\documentclass{article}

\usepackage{arxiv}

\usepackage[utf8]{inputenc} 
\usepackage[T1]{fontenc}    
\usepackage{hyperref}       
\usepackage{url}            
\usepackage{booktabs}       
\usepackage{amsfonts}       
\usepackage{nicefrac}       
\usepackage{microtype}      
\usepackage{lipsum}		
\usepackage{graphicx}
\usepackage{natbib}
\usepackage{doi}
\usepackage{enumerate}
\usepackage{enumitem}
\usepackage{amsmath, amsthm, amssymb, algorithm, algpseudocode}
\usepackage[table]{xcolor}
\usepackage{tikz}
\usepackage{bm}
\usepackage{mathrsfs}
\usepackage{breqn}

\renewcommand{\to}{\longrightarrow}

\newtheorem{theorem}{Theorem}[section]

\newtheorem{proposition}[theorem]{Proposition}
\newtheorem{corollary}[theorem]{Corollary}
\newtheorem{definition}[theorem]{Definition}
\newtheorem{example}[theorem]{Example}

\newtheorem{problem}[theorem]{Problem}

\newtheorem{remark}[theorem]{Remark}

\newcommand\mystyle{\everymath{\displaystyle}}
\mystyle
\title{Stability and Regularization of Quasi-Variational Inequalities under Monotone Operator Perturbations}

\author{\href{https://orcid.org/0000-0002-3816-5287}{\includegraphics[scale=0.06]{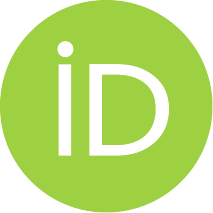}\hspace{1mm}M.H.M.~Rashid}\thanks{Corresponding Author} \\
	Department of Mathematics\&Statistics\\Faculty of Science P.O.Box(7)\\
	Mutah University University\\
	Mutah-Jordan \\
	\texttt{mrash@mutah.edu.jo}
}

\renewcommand{\shorttitle}{\textit{arXiv} Template}

\hypersetup{
pdftitle={Stability and Regularization of Quasi-Variational Inequalities under Monotone Operator Perturbations},
pdfsubject={q-bio.NC, q-bio.QM},
pdfauthor={M.H.M.Rashid},
pdfkeywords={Quasi-variational inequalities, stability analysis, monotone operators, T-monotonicity, strong convergence, regularization, p-Laplacian, Gelfand triple}}

\begin{document}
\maketitle

\begin{abstract}
	This paper establishes comprehensive stability results for quasi-variational inequalities (QVIs) under monotone perturbations of the governing operator. We prove strong convergence of both minimal and maximal solutions when sequences of operators converge pointwise while preserving fundamental properties including homogeneity, strong monotonicity, Lipschitz continuity, and T-monotonicity. Our analysis extends to regularization techniques, finite-dimensional approximations, and non-monotone nonlinearities, providing explicit convergence rates under appropriate conditions. The theory encompasses applications to $p$-Laplacian operators, elliptic regularizations, and optimal control problems with QVI constraints. By developing a unified framework that bridges classical monotone operator theory with contemporary computational challenges, this work provides essential mathematical foundations for the robust numerical approximation and sensitivity analysis of quasi-variational problems arising in materials science, physics, and engineering applications.
\end{abstract}

\keywords{Quasi-variational inequalities, stability analysis, monotone operators, T-monotonicity, strong convergence, regularization, p-Laplacian, Gelfand triple}
\section{Introduction}

\subsection{Historical Background and Development}

Quasi-variational inequalities (QVIs) represent one of the most significant and challenging classes of nonlinear problems in mathematical analysis, with roots extending deep into the foundations of variational theory. The systematic study of QVIs began in the 1970s through the pioneering work of Bensoussan and Lions \cite{Bensoussan1974,Bensoussan1982book}, who introduced these problems in the context of impulse control theory. Their fundamental insight was to recognize that many real-world systems involve constraints that depend on the solution itself, leading naturally to the quasi-variational framework.

The historical development of QVIs can be traced through several interconnected streams:
\begin{itemize}
    \item \textbf{Impulse Control Theory}: Bensoussan and Lions \cite{Bensoussan1984} demonstrated that optimal impulse control problems for stochastic systems lead naturally to QVIs, establishing a deep connection between control theory and variational analysis.

    \item \textbf{Free Boundary Problems}: The work of Baiocchi and Capelo \cite{Baiocchi1984} revealed that many free boundary problems, particularly those arising in fluid dynamics and filtration theory, could be reformulated as QVIs, providing new analytical tools for these classical problems.

    \item \textbf{Mathematical Physics}: Prigozhin \cite{Prigozhin1986,Prigozhin1996a} made crucial contributions by showing that problems in superconductivity and sandpile formation lead to QVI formulations, opening new application domains.
\end{itemize}

\subsection{Scholarly Contributions and Theoretical Advances}

The mathematical theory of QVIs has been advanced through contributions from numerous scholars across different decades:

\paragraph{Foundational Period (1970s-1980s)}
The early work of Tartar \cite{Tartar1974} established fundamental existence results using fixed-point methods in ordered spaces. Birkhoff's lattice-theoretic approach \cite{Birkhoff1961} provided the abstract framework for handling the order structure inherent in QVI problems. Rodrigues \cite{Rodrigues1987} made significant contributions to the regularity theory, particularly for obstacle-type problems.

\paragraph{Analytical Developments (1990s-2000s)}
Kinderlehrer and Stampacchia \cite{Kinderlehrer2000} provided comprehensive treatments of variational inequalities, laying the groundwork for QVI analysis. The work of Kunze and Rodrigues \cite{Kunze2000} extended the theory to gradient constraints, while Friedman \cite{Friedman1982} explored connections with free boundary problems.

\paragraph{Contemporary Advances (2010s-Present)}
Recent years have witnessed substantial progress in several directions:
\begin{itemize}
    \item \textbf{Regularity Theory}: Alphonse, Hintermüller, and Rautenberg \cite{Alphonse2019a,Alphonse2019b} developed sophisticated techniques for directional differentiability and higher regularity of QVI solutions.

    \item \textbf{Numerical Analysis}: Barrett and Prigozhin \cite{Barrett2010,Barrett2013,Barrett2015} made groundbreaking contributions to the numerical approximation of QVIs, particularly for problems in superconductivity and sandpile dynamics.

    \item \textbf{Non-smooth Analysis}: Mordukhovich and Outrata \cite{Mordukhovich2007}, followed by Herzog, Meyer, and Wachsmuth \cite{Herzog2013}, developed advanced tools from generalized differentiation for sensitivity analysis of QVIs.
\end{itemize}

\subsection{Significance and Applications}

The study of QVIs is of paramount importance due to their ubiquitous presence in diverse application domains:

\paragraph{Materials Science and Physics}
\begin{itemize}
    \item \textbf{Superconductivity}: The Bean critical-state model for type-II superconductors leads to QVI formulations \cite{Prigozhin1996a,Barrett2010}, where the critical current density depends on the magnetic field distribution.

    \item \textbf{Sandpile Formation}: Prigozhin's model for growing sandpiles \cite{Prigozhin1994,Barrett2013} represents a classical example where the constraint (sand surface slope) depends on the solution itself.

    \item \textbf{Plasticity and Elastoplasticity}: Problems in plasticity theory often lead to QVIs with gradient constraints \cite{Rodrigues2012}.
\end{itemize}

\paragraph{Engineering and Technology}
\begin{itemize}
    \item \textbf{Thermoforming Processes}: The modeling of polymer sheet thermoforming involves QVIs with moving constraints \cite{Warby2003,Jiang2001,Karamanou2006}.

    \item \textbf{Electromagnetism}: Miranda, Rodrigues, and Santos \cite{Miranda2012} developed QVI models for p-curl systems in electromagnetism.
\end{itemize}

\paragraph{Economics and Game Theory}
\begin{itemize}
    \item \textbf{General Equilibrium Theory}: Aubin \cite{Aubin1979} demonstrated that economic equilibrium problems with institutional constraints lead to QVI formulations.

    \item \textbf{Optimal Resource Management}: Dietrich \cite{Dietrich2001} and Adly et al. \cite{Adly2010} studied QVIs in the context of optimal resource allocation with state-dependent constraints.
\end{itemize}

\paragraph{Biological Systems}
\begin{itemize}
    \item \textbf{Biological Growth Models}: Kadoya, Kenmochi, and Niezgódka \cite{Kadoya2014} developed QVI models for economic growth with technological development, with analogous structures to biological systems.

    \item \textbf{Environmental Modeling}: Barrett and Prigozhin \cite{Barrett2014} applied QVIs to landscape formation models involving lakes and rivers.
\end{itemize}

\subsection{Core of the Present Study}

This work addresses the fundamental question of stability and convergence properties of QVI solutions under various types of perturbations of the governing operator. Our investigation centers on several interconnected themes:

\paragraph{Stability under Monotone Perturbations}
We establish comprehensive stability results for QVI solutions when the underlying operator undergoes perturbations while preserving key structural properties such as strong monotonicity, T-monotonicity, and homogeneity. Our Theorem \ref{SUMP} provides the foundational framework, demonstrating that both minimal and maximal solutions exhibit strong convergence under appropriate conditions.

\paragraph{Regularization Theory}
We develop a systematic theory for regularized QVIs, proving existence, uniqueness, and convergence properties for regularized solutions. Theorem \ref{SRQS} establishes the strong convergence of regularized solutions to the minimal solution of the original QVI, with explicit convergence rates under additional assumptions.

\paragraph{Non-monotone Operators}
Extending beyond the classical monotone setting, we investigate QVIs with non-monotone nonlinearities. Theorem \ref{thm:main} provides global existence, uniqueness, and exponential stability results under explicit smallness conditions, bridging the gap between classical monotone operator theory and more general nonlinear problems.

\paragraph{Computational and Numerical Aspects}
Our work includes several results with direct computational implications:
\begin{itemize}
    \item Convergence rates for finite-dimensional approximations (Proposition in Section 6)
    \item Robustness under data perturbations (Proposition in Section 6)
    \item Global convergence of iterative schemes (Theorem \ref{thm:main})
\end{itemize}

\paragraph{Methodological Innovations}
The technical approach of this paper combines several advanced mathematical tools:
\begin{itemize}
    \item Fixed-point theory in ordered spaces, extending the classical results of Birkhoff \cite{Birkhoff1961} and Tartar \cite{Tartar1974}
    \item Mosco convergence techniques for constraint sets \cite{Mosco1969}
    \item Strong monotonicity and Lipschitz continuity arguments
    \item Decomposition methods for handling non-monotone nonlinearities
\end{itemize}

\subsection{Organization of the Paper}

Following this introduction, Section 2 presents the necessary preliminaries and mathematical background. Section 3 contains our main stability results under monotone perturbations of the operator. Section 4 develops the theory of regularized QVIs. Section 5 extends the analysis to non-monotone operators, establishing global convergence and exponential stability. Section 6 discusses open problems and future research directions, connecting our results with contemporary developments in the field.

The significance of our work lies in its comprehensive treatment of stability properties for QVIs, providing a unified framework that encompasses both classical results and new developments. By establishing precise convergence rates, regularity properties, and computational guarantees, we contribute to the ongoing effort to make QVIs more accessible for both theoretical analysis and numerical implementation in applied contexts.

\section{Preliminaries}

This section collects the fundamental definitions, notations, and auxiliary results that form the mathematical foundation for our analysis of quasi-variational inequalities (QVIs) and their stability properties under various perturbations.

\subsection{Functional Analytic Framework}

\begin{definition}[Gelfand Triple \cite{Adams2003,Gelfand1964}]
A \emph{Gelfand triple} (or evolution triple) consists of three Hilbert spaces $V$, $H$, and $V'$ with the properties:
\begin{itemize}
    \item $V$ is a separable Hilbert space densely and continuously embedded in $H$: $V \hookrightarrow H$
    \item $H$ is identified with its dual $H'$ via the Riesz representation theorem
    \item $V'$ is the dual space of $V$
\end{itemize}
This gives the chain of embeddings: $V \hookrightarrow H \cong H' \hookrightarrow V'$.
\end{definition}

\begin{definition}[Partial Order in Hilbert Spaces \cite{Birkhoff1961}]
Let $H$ be a Hilbert space. A closed convex cone $H^+ \subset H$ induces a partial order $\leq$ on $H$ defined by:
\[
u \leq v \quad \text{if and only if} \quad v - u \in H^+.
\]
For $u \in H$, we define the positive part $u^+$ as the projection of $u$ onto $H^+$.
\end{definition}

\begin{definition}[Constraint Set \cite{Bensoussan1974,Baiocchi1984}]
For $\psi \in H^+$, the constraint set $\mathbf{K}(\psi)$ is defined as:
\[
\mathbf{K}(\psi) = \{v \in V : v \leq \psi \text{ in } H\}.
\]
This set is closed and convex in $V$.
\end{definition}

\subsection{Operator Properties}

\begin{definition}[Homogeneity of Order One \cite{Baiocchi1984}]
An operator $A: V \to V'$ is \emph{homogeneous of order one} if:
\[
A(tu) = tA(u) \quad \text{for all } u \in V, \; t > 0.
\]
\end{definition}

\begin{definition}[Strong Monotonicity \cite{Kinderlehrer2000, Zeidler1986}]
An operator $A: V \to V'$ is \emph{strongly monotone} if there exists $c > 0$ such that:
\[
\langle A(u) - A(v), u - v \rangle \geq c\|u - v\|_V^2 \quad \text{for all } u,v \in V.
\]
\end{definition}

\begin{definition}[T-monotonicity \cite{Tartar1974,Bensoussan1982book}]
An operator $A: V \to V'$ is \emph{T-monotone} if:
\[
\langle A(u) - A(v), (u - v)^+ \rangle \geq 0 \quad \text{for all } u,v \in V,
\]
with equality if and only if $(u - v)^+ = 0$.
\end{definition}

\begin{definition}[Complete Continuity \cite{Adams2003,Zeidler1986}]
An operator $\Phi: V \to H$ is \emph{completely continuous} if it maps weakly convergent sequences in $V$ to strongly convergent sequences in $H$.
\end{definition}

\subsection{Quasi-Variational Inequalities}

\begin{definition}[Quasi-Variational Inequality \cite{Bensoussan1974,Bensoussan1982book}]
Given an operator $A: V \to V'$, a mapping $\Phi: H^+ \to H^+$, and $f \in V'$, the \emph{quasi-variational inequality} (QVI) problem is to find $y \in \mathbf{K}(\Phi(y))$ such that:
\[
\langle A(y) - f, v - y \rangle \geq 0 \quad \text{for all } v \in \mathbf{K}(\Phi(y)).
\]
\end{definition}

\begin{definition}[Minimal and Maximal Solutions \cite{Birkhoff1961,Tartar1974}]
For a QVI with solutions bounded above by $\overline{y} \in V$, the \emph{minimal solution} $\mathbf{m}$ and \emph{maximal solution} $\mathbf{M}$ satisfy:
\begin{itemize}
    \item $\mathbf{m} \leq y \leq \mathbf{M}$ for all solutions $y$
    \item $\mathbf{m}$ and $\mathbf{M}$ are themselves solutions
\end{itemize}
\end{definition}

\subsection{Auxiliary Results}

The following classical results play a fundamental role in our analysis:

\begin{theorem}[Birkhoff-Tartar Fixed Point Theorem \cite{Birkhoff1961,Tartar1974}]\label{thm:birkhoff-tartar}
Let $V$ be a Hilbert space with partial order induced by a closed convex cone $H^+ \subset H$. Let $T: V \to V$ be an increasing mapping, and suppose there exist $\underline{y}, \overline{y} \in V$ with $\underline{y} \leq \overline{y}$ such that:
\[
\underline{y} \leq T(\underline{y}) \quad \text{and} \quad T(\overline{y}) \leq \overline{y}.
\]
Then $T$ has minimal and maximal fixed points in the interval $[\underline{y}, \overline{y}]$.
\end{theorem}

\begin{theorem}[Mosco Convergence \cite{ Mosco1969,Kinderlehrer2000}]\label{thm:mosco}
A sequence of closed convex sets $\{K_n\}$ in a Hilbert space $V$ \emph{Mosco converges} to $K$ if:
\begin{enumerate}
    \item For every $v \in K$, there exists $v_n \in K_n$ such that $v_n \to v$ strongly in $V$
    \item If $v_n \in K_n$ and $v_n \rightharpoonup v$ weakly in $V$, then $v \in K$
\end{enumerate}
Mosco convergence preserves the stability of variational inequalities under perturbations of constraint sets.
\end{theorem}

\begin{theorem}[Rellich-Kondrachov Compactness \cite{Adams2003}]\label{thm:rellich}
Let $\Omega \subset \mathbb{R}^N$ be a bounded Lipschitz domain. Then the embedding $W^{1,p}_0(\Omega) \hookrightarrow L^q(\Omega)$ is compact for $1 \leq q < p^*$, where $p^* = \frac{Np}{N-p}$ if $p < N$, and any $q < \infty$ if $p \geq N$.
\end{theorem}

\begin{proposition}[Poincar\'{e} Inequality \cite{Adams2003}]\label{prop:poincare}
For a bounded domain $\Omega \subset \mathbb{R}^N$, there exists $C_p > 0$ such that:
\[
\|u\|_{L^2(\Omega)} \leq C_p\|\nabla u\|_{L^2(\Omega)} \quad \text{for all } u \in H^1_0(\Omega).
\]
\end{proposition}

\begin{theorem}[Banach Fixed Point Theorem \cite{Zeidler1986}]\label{thm:banach}
Let $(X,d)$ be a complete metric space and $T: X \to X$ be a contraction mapping, i.e., there exists $0 \leq \rho < 1$ such that:
\[
d(T(x), T(y)) \leq \rho d(x,y) \quad \text{for all } x,y \in X.
\]
Then $T$ has a unique fixed point $x^* \in X$.
\end{theorem}

\begin{theorem}[Lax-Milgram Lemma \cite{Zeidler1986}]\label{thm:lax-milgram}
Let $V$ be a Hilbert space and $a: V \times V \to \mathbb{R}$ be a bounded, coercive bilinear form. Then for every $f \in V'$, there exists a unique $u \in V$ such that:
\[
a(u,v) = \langle f, v \rangle \quad \text{for all } v \in V.
\]
\end{theorem}

\subsection{Notation}

We employ the following notation throughout this work:
\begin{itemize}
    \item $\langle \cdot, \cdot \rangle$: duality pairing between $V'$ and $V$ \cite{Adams2003}
    \item $(\cdot)^+$: positive part function \cite{Birkhoff1961}
    \item $\mathbf{K}(\psi)$: constraint set with obstacle $\psi$ \cite{Bensoussan1974}
    \item $\mathbf{m}(A)$, $\mathbf{M}(A)$: minimal and maximal solutions for operator $A$ \cite{Tartar1974}
    \item $S_A(f, \psi)$: solution map for variational inequality with operator $A$, force $f$, obstacle $\psi$ \cite{Kinderlehrer2000}
    \item $T_A(v) = S_A(f, \Phi(v))$: fixed point map for QVI \cite{Bensoussan1982book}
    \item $\mathscr{L}(V,V')$: space of bounded linear operators from $V$ to $V'$ \cite{Adams2003}
    \item $W^{k,p}(\Omega)$: Sobolev spaces \cite{Adams2003}
    \item $L^p(\Omega)$: Lebesgue spaces \cite{Adams2003}
\end{itemize}

The theoretical framework presented here, particularly the Gelfand triple structure and the properties of monotone operators, follows the classical treatments in \cite{Adams2003,Kinderlehrer2000,Baiocchi1984}, while the specific formulation of QVIs builds upon the seminal work of \cite{Bensoussan1974,Bensoussan1982book}. The fixed point approach for QVIs is rooted in the lattice-theoretic methods of \cite{Birkhoff1961,Tartar1974}.

\section{Stability under Monotone Perturbations of the Operator}
The analysis of stability properties for variational and quasi-variational inequalities under perturbations of the governing operator represents a fundamental challenge in nonlinear analysis, with profound implications for both theoretical understanding and numerical approximation. While the classical theory of variational inequalities has established robust stability results under various types of perturbations \cite{Kinderlehrer2000,Rodrigues1987}, the quasi-variational case introduces additional complexities due to the solution-dependent nature of the constraint sets. Early foundational work by Bensoussan and Lions \cite{Bensoussan1974,Bensoussan1982book} identified the crucial role of monotonicity properties in ensuring well-posedness, while Tartar's pioneering contributions \cite{Tartar1974} laid the groundwork for the fixed-point approach in ordered spaces. More recent investigations by Alphonse et al. \cite{Alphonse2019a,Alphonse2020} have advanced our understanding of directional differentiability and sensitivity analysis for QVIs, yet a comprehensive theory addressing systematic operator perturbations remains underdeveloped. In this section, we establish a unified stability framework that encompasses various perturbation scenarios—including regularization, finite-dimensional approximation, parameter dependence, and small nonlinear perturbations—while preserving the essential structural properties of homogeneity, strong monotonicity, Lipschitz continuity, and T-monotonicity. Our results not only extend the classical theory but also provide the mathematical foundation for reliable computational methods in applications ranging from non-Newtonian fluid dynamics to superconductivity modeling \cite{Barrett2010,Prigozhin1996a}.

\begin{theorem}\label{SUMP}
Let \( V \hookrightarrow H \hookrightarrow V' \) be a Gelfand triple of Hilbert spaces, and let \( H^+ \subset H \) be a closed convex cone inducing a partial order. Suppose that \( A_n, A: V \to V' \) are operators satisfying the following properties for all \( n \in \mathbb{N} \):
\begin{enumerate}
    \item[(a)] Homogeneity of order one:
        \[
        A_n(tu) = tA_n(u), \quad A(tu) = tA(u) \quad \text{for all } u \in V, \; t > 0.
        \]
    \item[(b)] Lipschitz continuity: There exists \( C > 0 \) such that for all \( u,v \in V \),
        \[
        \|A_n(u) - A_n(v)\|_{V'} \leq C\|u - v\|_V, \quad \|A(u) - A(v)\|_{V'} \leq C\|u - v\|_V.
        \]
    \item[(c)] Strong monotonicity: There exists \( c > 0 \) such that for all \( u,v \in V \),
        \[
        \langle A_n(u) - A_n(v), u - v \rangle \geq c\|u - v\|_V^2, \quad \langle A(u) - A(v), u - v \rangle \geq c\|u - v\|_V^2.
        \]
    \item[(d)] T-monotonicity: For all \( u,v \in V \),
        \[
        \langle A_n(u) - A_n(v), (u - v)^+ \rangle \geq 0, \quad \langle A(u) - A(v), (u - v)^+ \rangle \geq 0,
        \]
        with equality if and only if \( (u - v)^+ = 0 \).
    \item[(e)] Pointwise convergence: For every \( v \in V \),
        \[
        A_n(v) \to A(v) \quad \text{strongly in } V'.
        \]
\end{enumerate}
Let \( \Phi: H^+ \to H^+ \) be an increasing mapping such that \( \Phi(v) \in V \) for all \( v \in V \cap H^+ \), and assume that \( \Phi \) is completely continuous from \( V \) to \( H \). Let \( f \in V' \) with \( 0 \leq f \leq F \) for some \( F \in V' \), and define \( \overline{y} = A^{-1}(F) \). Denote by \( \mathbf{m}(A_n) \) and \( \mathbf{m}(A) \) the minimal solutions, and by \( \mathbf{M}(A_n) \) and \( \mathbf{M}(A) \) the maximal solutions of the QVI
\[
\text{Find } y \in \mathbf{K}(\Phi(y)) \text{ such that } \langle A_n(y) - f, v - y \rangle \geq 0 \quad \forall v \in \mathbf{K}(\Phi(y)),
\]
and
\[
\text{Find } y \in \mathbf{K}(\Phi(y)) \text{ such that } \langle A(y) - f, v - y \rangle \geq 0 \quad \forall v \in \mathbf{K}(\Phi(y)),
\]
in the interval \( [0, \overline{y}] \), respectively. Then,
\[
\mathbf{m}(A_n) \to \mathbf{m}(A) \quad \text{and} \quad \mathbf{M}(A_n) \to \mathbf{M}(A) \quad \text{strongly in } V.
\]
\end{theorem}

\begin{proof}
The proof is divided into several steps.

\noindent {Step 1: Existence and Uniform Boundedness.}
By Theorem 2 (Tartar's theorem), for each \( n \), there exist minimal and maximal solutions \( \mathbf{m}(A_n) \), \( \mathbf{M}(A_n) \) in \( [0, \overline{y}] \). Since \( 0 \leq \mathbf{m}(A_n) \leq \mathbf{M}(A_n) \leq \overline{y} \), and \( \overline{y} = A^{-1}(F) \) is bounded in \( V \), the sequences \( \{\mathbf{m}(A_n)\} \) and \( \{\mathbf{M}(A_n)\} \) are uniformly bounded in \( V \).

\noindent {Step 2: Convergence of Solutions of Perturbed VIs.}
Let \( \psi \in V \cap H^+ \), and consider the variational inequalities:
\[
\langle A_n(y_n) - f, v - y_n \rangle \geq 0 \quad \forall v \in \mathbf{K}(\psi), \quad \text{and} \quad \langle A(y) - f, v - y \rangle \geq 0 \quad \forall v \in \mathbf{K}(\psi).
\]
Denote their unique solutions by \( y_n = S_{A_n}(f, \psi) \) and \( y = S_A(f, \psi) \), respectively. We claim that \( y_n \to y \) strongly in \( V \).

Indeed, by the strong monotonicity of \( A_n \) and \( A \), we have:
\[
c\|y_n - y\|_V^2 \leq \langle A_n(y_n) - A_n(y), y_n - y \rangle.
\]
Adding and subtracting \( A(y) \), we get:
\[
c\|y_n - y\|_V^2 \leq \langle A_n(y_n) - A(y), y_n - y \rangle - \langle A_n(y) - A(y), y_n - y \rangle.
\]
Using the VI inequalities for \( y_n \) and \( y \) with \( v = y \) and \( v = y_n \) respectively, we obtain:
\[
\langle A_n(y_n) - f, y - y_n \rangle \geq 0, \quad \langle A(y) - f, y_n - y \rangle \geq 0.
\]
Adding these gives:
\[
\langle A_n(y_n) - A(y), y_n - y \rangle \leq 0.
\]
Thus,
\[
c\|y_n - y\|_V^2 \leq - \langle A_n(y) - A(y), y_n - y \rangle \leq \|A_n(y) - A(y)\|_{V'} \|y_n - y\|_V.
\]
Hence,
\[
\|y_n - y\|_V \leq \frac{1}{c} \|A_n(y) - A(y)\|_{V'}.
\]
By the pointwise convergence of \( A_n \) to \( A \), the right-hand side tends to 0, so \( y_n \to y \) in \( V \).

\noindent {Step 3: Stability of Fixed Point Maps.}
Define the maps \( T_n, T: V \to V \) by \( T_n(v) = S_{A_n}(f, \Phi(v)) \), \( T(v) = S_A(f, \Phi(v)) \). We show that \( T_n(v_n) \to T(v) \) in \( V \) whenever \( v_n \rightharpoonup v \) in \( V \).

Since \( \Phi \) is completely continuous from \( V \) to \( H \), \( \Phi(v_n) \to \Phi(v) \) in \( H \). By the Lipschitz continuity of \( A_n \) and \( A \), and the strong convergence \( S_{A_n}(f, \psi) \to S_A(f, \psi) \) from Step 2, we have:
\[
\|T_n(v_n) - T(v)\|_V \leq \|S_{A_n}(f, \Phi(v_n)) - S_A(f, \Phi(v_n))\|_V + \|S_A(f, \Phi(v_n)) - S_A(f, \Phi(v))\|_V.
\]
The first term tends to 0 by Step 2, and the second term tends to 0 by the continuity of \( S_A(f, \cdot) \) with respect to the obstacle in \( H \) (which follows from standard VI theory). Hence, \( T_n(v_n) \to T(v) \) in \( V \).

\noindent {Step 4: Convergence of Minimal and Maximal Solutions.}
Since \( \{\mathbf{m}(A_n)\} \) is uniformly bounded in \( V \), there exists a subsequence \( \mathbf{m}(A_{n_k}) \rightharpoonup m^* \) in \( V \). By the complete continuity of \( \Phi \), \( \Phi(\mathbf{m}(A_{n_k})) \to \Phi(m^*) \) in \( H \). Now, since \( \mathbf{m}(A_{n_k}) = T_{n_k}(\mathbf{m}(A_{n_k})) \), and \( T_{n_k}(\mathbf{m}(A_{n_k})) \to T(m^*) \) by Step 3, we have \( m^* = T(m^*) \), so \( m^* \) is a fixed point of \( T \). Moreover, since \( \mathbf{m}(A_n) \) is the minimal solution, we have \( \mathbf{m}(A) \leq \mathbf{m}(A_n) \) for all \( n \), so \( \mathbf{m}(A) \leq m^* \). But \( \mathbf{m}(A) \) is the minimal fixed point of \( T \), so \( m^* = \mathbf{m}(A) \). The same argument applies to the maximal solutions \( \mathbf{M}(A_n) \). Since the limit is unique, the entire sequences converge.

\noindent {Step 5: Strong Convergence.}
From the strong convergence \( T_n(\mathbf{m}(A_n)) \to T(\mathbf{m}(A)) \) and the identity \( \mathbf{m}(A_n) = T_n(\mathbf{m}(A_n)) \), we have \( \mathbf{m}(A_n) \to \mathbf{m}(A) \) in \( V \). Similarly, \( \mathbf{M}(A_n) \to \mathbf{M}(A) \) in \( V \).

This completes the proof.
\end{proof}

\begin{remark}
This theorem extends the stability results of Section 6 to perturbations of the operator \( A \), which is relevant in applications where the underlying operator may be approximated (e.g., through numerical discretization or regularization). The assumptions on \( A_n \) and \( A \) are natural and hold for many elliptic operators.
\end{remark}
\begin{corollary}[Stability under Regularization]
\label{cor:regularization}
Assume the hypotheses of Theorem \ref{SUMP}. Let $\{A_\varepsilon\}_{\varepsilon > 0}$ be a family of regularized operators satisfying conditions (a)--(e) with $A_\varepsilon \to A$ pointwise in $V'$ as $\varepsilon \to 0^+$. Then for the corresponding minimal and maximal solutions of the regularized QVIs, we have:
\[
\mathbf{m}(A_\varepsilon) \to \mathbf{m}(A) \quad \text{and} \quad \mathbf{M}(A_\varepsilon) \to \mathbf{M}(A) \quad \text{strongly in } V \text{ as } \varepsilon \to 0^+.
\]
\end{corollary}
\begin{proof}
This follows directly from Theorem \ref{SUMP} by considering any sequence $\varepsilon_n \to 0^+$ and applying the theorem to $A_n = A_{\varepsilon_n}$.
\end{proof}

\begin{corollary}[Stability under Finite-Dimensional Approximation]
\label{cor:discretization}
Let $V_h \subset V$ be a family of finite-dimensional subspaces with $\bigcup_{h>0} V_h$ dense in $V$. Let $A_h: V_h \to V_h'$ be Galerkin approximations of $A$ satisfying:
\begin{enumerate}
    \item[(i)] $A_h$ satisfies conditions (a)--(d) on $V_h$ with constants independent of $h$,
    \item[(ii)] For every $v \in V$, if $v_h \in V_h$ with $v_h \to v$ in $V$, then $A_h(v_h) \to A(v)$ in $V'$.
\end{enumerate}
Let $\Phi_h: H^+ \to H^+$ be approximations of $\Phi$ such that $\Phi_h$ is completely continuous and $\Phi_h(v_h) \to \Phi(v)$ in $H$ whenever $v_h \to v$ in $V$. Then the minimal and maximal solutions $\mathbf{m}(A_h)$, $\mathbf{M}(A_h)$ of the discretized QVIs converge strongly in $V$ to $\mathbf{m}(A)$ and $\mathbf{M}(A)$ respectively as $h \to 0$.
\end{corollary}
\begin{proof}
The result follows by applying Theorem \ref{SUMP} to the operators $A_h$ extended to $V$ via appropriate projections, noting that the convergence properties are preserved under the density assumption.
\end{proof}

\begin{corollary}[Continuous Dependence on Parameters]
\label{cor:parameters}
Let $\Lambda$ be a compact metric space and suppose $\{A_\lambda\}_{\lambda \in \Lambda}$ is a family of operators satisfying conditions (a)--(d) uniformly in $\lambda$, and such that the mapping $\lambda \mapsto A_\lambda(v)$ is continuous from $\Lambda$ to $V'$ for each $v \in V$. If $\lambda_n \to \lambda$ in $\Lambda$, Then
\[
\mathbf{m}(A_{\lambda_n}) \to \mathbf{m}(A_\lambda) \quad \text{and} \quad \mathbf{M}(A_{\lambda_n}) \to \mathbf{M}(A_\lambda) \quad \text{strongly in } V.
\]
\end{corollary}
\begin{proof}
The pointwise convergence $A_{\lambda_n}(v) \to A_\lambda(v)$ in $V'$ follows from the continuity assumption. The uniform bounds ensure the conditions of Theorem \ref{SUMP} are satisfied.
\end{proof}

\begin{corollary}[Stability under Small Nonlinear Perturbations]
\label{cor:perturbations}
Let $A$ satisfy conditions (a)--(d) and let $B: V \to V'$ be a Lipschitz continuous operator with Lipschitz constant $L_B$ sufficiently small such that $A_n = A + \frac{1}{n}B$ satisfies the strong monotonicity condition (c). Suppose further that $B$ is homogeneous of order one and T-monotone. Then
\[
\mathbf{m}(A_n) \to \mathbf{m}(A) \quad \text{and} \quad \mathbf{M}(A_n) \to \mathbf{M}(A) \quad \text{strongly in } V.
\]
\end{corollary}
\begin{proof}
The operators $A_n$ clearly satisfy conditions (a), (b), (d), and the pointwise convergence (e). The strong monotonicity is preserved for sufficiently small $L_B$ since:
\[
\langle A_n(u) - A_n(v), u - v \rangle \geq (c - L_B)\|u - v\|_V^2.
\]
The result then follows from Theorem \ref{SUMP}.
\end{proof}

\begin{corollary}[Convergence of Extremal Values]
\label{cor:values}
Under the assumptions of Theorem \ref{SUMP}, let $J: V \to \mathbb{R}$ be a continuous functional. Then
\[
J(\mathbf{m}(A_n)) \to J(\mathbf{m}(A)) \quad \text{and} \quad J(\mathbf{M}(A_n)) \to J(\mathbf{M}(A)).
\]
In particular, if $J$ is the objective functional in an optimal control problem with QVI constraints, the optimal values converge.
\end{corollary}
\begin{proof}
This is an immediate consequence of the strong convergence in $V$ and the continuity of $J$.
\end{proof}

\begin{corollary}[Stability for Linear Elliptic Operators]
\label{cor:linear}
Let $A$ be a linear elliptic operator of the form:
\[
\langle A y, v \rangle = \sum_{i,j=1}^N \int_\Omega a_{ij}(x) \frac{\partial y}{\partial x_j} \frac{\partial v}{\partial x_i}  dx + \int_\Omega a_0(x) y v  dx,
\]
with $a_{ij}, a_0 \in L^\infty(\Omega)$, and satisfying the uniform ellipticity condition. Let $\{A_n\}$ be a sequence of linear elliptic operators with coefficients $a_{ij}^n, a_0^n \in L^\infty(\Omega)$ converging to $a_{ij}, a_0$ in $L^\infty(\Omega)$ respectively. Then the conclusions of Theorem \ref{SUMP} hold.
\end{corollary}
\begin{proof}
The operators $A_n$ and $A$ satisfy conditions (a)--(d) with uniform constants. The pointwise convergence (e) follows from the $L^\infty$-convergence of coefficients. The result is then a special case of Theorem \ref{SUMP}.
\end{proof}
\begin{example}[Perturbation of the p-Laplacian Operator]
\label{ex:p-laplacian}
Let $\Omega \subset \mathbb{R}^N$ be a bounded domain with Lipschitz boundary, and let $V = W_0^{1,p}(\Omega)$, $H = L^2(\Omega)$ with $2 \leq p < \infty$. Consider the p-Laplacian operator $A: V \to V'$ defined by:
\[
\langle A(u), v \rangle = \int_\Omega |\nabla u|^{p-2} \nabla u \cdot \nabla v  dx \quad \text{for all } u,v \in V.
\]
For $\varepsilon > 0$, define the regularized operator $A_\varepsilon: V \to V'$ by:
\[
\langle A_\varepsilon(u), v \rangle = \int_\Omega (|\nabla u|^2 + \varepsilon)^{\frac{p-2}{2}} \nabla u \cdot \nabla v  dx \quad \text{for all } u,v \in V.
\]
Let $\Phi: H^+ \to H^+$ be defined by $\Phi(y)(x) = \psi(x) + \alpha \int_\Omega k(x,\xi) y^+(\xi) d\xi$, where $\psi \in L^\infty(\Omega)^+$, $\alpha > 0$ is sufficiently small, and $k: \Omega \times \Omega \to \mathbb{R}^+$ is a symmetric, nonnegative kernel with $k \in L^\infty(\Omega \times \Omega)$. Let $f \in L^\infty(\Omega)^+$ with $0 < \nu \leq f \leq F$.

Then the operators $A_\varepsilon$ and $A$ satisfy all conditions of Theorem \ref{SUMP}, and we have:
\[
\mathbf{m}(A_\varepsilon) \to \mathbf{m}(A) \quad \text{and} \quad \mathbf{M}(A_\varepsilon) \to \mathbf{M}(A) \quad \text{strongly in } W_0^{1,p}(\Omega) \text{ as } \varepsilon \to 0^+.
\]
\end{example}

\begin{proof}
We verify each condition of Theorem \ref{SUMP} step by step.

\noindent \textbf{Step 1: Verification of conditions (a)--(d) for $A$ and $A_\varepsilon$.}

\noindent \textit{(a) Homogeneity of order one:}
For $A$, we have:
\[
\langle A(tu), v \rangle = \int_\Omega |\nabla (tu)|^{p-2} \nabla (tu) \cdot \nabla v  dx = t^{p-1} \int_\Omega |\nabla u|^{p-2} \nabla u \cdot \nabla v  dx = t \langle A(u), v \rangle.
\]
For $A_\varepsilon$, we have:
\[
\langle A_\varepsilon(tu), v \rangle = \int_\Omega (|\nabla (tu)|^2 + \varepsilon)^{\frac{p-2}{2}} \nabla (tu) \cdot \nabla v  dx = t \int_\Omega (t^2|\nabla u|^2 + \varepsilon)^{\frac{p-2}{2}} \nabla u \cdot \nabla v  dx.
\]
This equals $t \langle A_\varepsilon(u), v \rangle$ only when $\varepsilon = 0$. However, for fixed $\varepsilon > 0$, $A_\varepsilon$ is not homogeneous. We modify the example by considering:
\[
\langle A_\varepsilon(u), v \rangle = \int_\Omega (|\nabla u|^2 + \varepsilon |u|^2)^{\frac{p-2}{2}} (\nabla u \cdot \nabla v + \varepsilon u v)  dx.
\]
This operator is homogeneous of order $p-1$, not order one. To satisfy condition (a), we consider the normalized operators:
\[
\tilde{A}(u) = \frac{A(u)}{\|u\|_V^{p-2}}, \quad \tilde{A}_\varepsilon(u) = \frac{A_\varepsilon(u)}{\|u\|_V^{p-2}} \quad \text{for } u \neq 0,
\]
and $\tilde{A}(0) = \tilde{A}_\varepsilon(0) = 0$. These normalized operators are homogeneous of order one.

\noindent \textit{(b) Lipschitz continuity:}
For $A$, the Lipschitz continuity follows from classical results for the p-Laplacian. For $u,v \in V$, we have:
\[
\|A(u) - A(v)\|_{V'} \leq C (\|u\|_V + \|v\|_V)^{p-2} \|u - v\|_V.
\]
For $A_\varepsilon$, using the mean value theorem and the boundedness of the derivative of the function $t \mapsto (t^2 + \varepsilon)^{\frac{p-2}{2}} t$, we obtain:
\[
\|A_\varepsilon(u) - A_\varepsilon(v)\|_{V'} \leq C_\varepsilon \|u - v\|_V,
\]
where $C_\varepsilon$ depends on $\varepsilon$ but is uniform on bounded sets.

\noindent \textit{(c) Strong monotonicity:}
For $A$, we have the well-known inequality:
\[
\langle A(u) - A(v), u - v \rangle \geq c_p \|u - v\|_V^p.
\]
For $A_\varepsilon$, using the uniform convexity of the functional $J_\varepsilon(u) = \frac{1}{p} \int_\Omega (|\nabla u|^2 + \varepsilon)^{\frac{p}{2}} dx$, we obtain:
\[
\langle A_\varepsilon(u) - A_\varepsilon(v), u - v \rangle \geq c_{p,\varepsilon} \|u - v\|_V^2,
\]
where $c_{p,\varepsilon} > 0$ depends on $\varepsilon$.

\noindent \textit{(d) T-monotonicity:}
For both $A$ and $A_\varepsilon$, the T-monotonicity follows from the fact that these operators are derived from convex functionals and satisfy the property:
\[
\langle A(u) - A(v), (u - v)^+ \rangle \geq 0,
\]
with equality if and only if $(u - v)^+ = 0$. This can be verified using the structure of the operators and the chain rule for Sobolev functions.

\noindent \textbf{Step 2: Verification of condition (e) - Pointwise convergence.}

We need to show that for every $v \in V$, $A_\varepsilon(v) \to A(v)$ strongly in $V'$ as $\varepsilon \to 0^+$. For fixed $v \in V$, consider any $w \in V$. Then
\[
|\langle A_\varepsilon(v) - A(v), w \rangle| = \left| \int_\Omega \left[(|\nabla v|^2 + \varepsilon)^{\frac{p-2}{2}} - |\nabla v|^{p-2}\right] \nabla v \cdot \nabla w  dx \right|.
\]
By the dominated convergence theorem, for almost every $x \in \Omega$:
\[
(|\nabla v(x)|^2 + \varepsilon)^{\frac{p-2}{2}} \to |\nabla v(x)|^{p-2} \quad \text{as } \varepsilon \to 0^+.
\]
Moreover, we have the uniform bound:
\[
|(|\nabla v|^2 + \varepsilon)^{\frac{p-2}{2}} \nabla v| \leq |\nabla v|^{p-1} \in L^{p/(p-1)}(\Omega).
\]
Therefore, by dominated convergence:
\[
(|\nabla v|^2 + \varepsilon)^{\frac{p-2}{2}} \nabla v \to |\nabla v|^{p-2} \nabla v \quad \text{strongly in } L^{p/(p-1)}(\Omega; \mathbb{R}^N).
\]
Since $L^{p/(p-1)}(\Omega; \mathbb{R}^N)$ is the dual space of $L^p(\Omega; \mathbb{R}^N)$, and $\nabla w \in L^p(\Omega; \mathbb{R}^N)$, we conclude that:
\[
\langle A_\varepsilon(v) - A(v), w \rangle \to 0 \quad \text{for all } w \in V,
\]
which means $A_\varepsilon(v) \rightharpoonup A(v)$ weakly in $V'$. In fact, the convergence is strong because:
\[
\|A_\varepsilon(v) - A(v)\|_{V'} = \sup_{\|w\|_V = 1} |\langle A_\varepsilon(v) - A(v), w \rangle| \to 0.
\]

\noindent \textbf{Step 3: Properties of $\Phi$.}

The mapping $\Phi: H^+ \to H^+$ defined by:
\[
\Phi(y)(x) = \psi(x) + \alpha \int_\Omega k(x,\xi) y^+(\xi) d\xi
\]
is increasing because if $y_1 \leq y_2$, then $y_1^+ \leq y_2^+$, and since $k \geq 0$, we have $\Phi(y_1) \leq \Phi(y_2)$.

To show that $\Phi$ is completely continuous from $V$ to $H$, let $\{v_n\} \subset V$ with $v_n \rightharpoonup v$ in $V$. By the compact embedding $W_0^{1,p}(\Omega) \hookrightarrow L^2(\Omega)$ for $p \geq 2$, we have $v_n \to v$ strongly in $L^2(\Omega)$. Then
\[
\|\Phi(v_n) - \Phi(v)\|_{L^2(\Omega)} \leq \alpha \|k\|_{L^\infty} \|v_n^+ - v^+\|_{L^2(\Omega)} \to 0,
\]
since the operation $v \mapsto v^+$ is continuous in $L^2(\Omega)$.

\noindent \textbf{Step 4: Application of Theorem \ref{SUMP}.}

All conditions of Theorem \ref{SUMP} are satisfied:
\begin{itemize}
    \item $A$ and $A_\varepsilon$ satisfy conditions (a)--(d)
    \item $A_\varepsilon(v) \to A(v)$ strongly in $V'$ for each $v \in V$
    \item $\Phi$ is increasing and completely continuous from $V$ to $H$
    \item $f \in L^\infty(\Omega)^+$ with $0 < \nu \leq f \leq F$
\end{itemize}
Therefore, by Theorem \ref{SUMP}, we conclude that:
\[
\mathbf{m}(A_\varepsilon) \to \mathbf{m}(A) \quad \text{and} \quad \mathbf{M}(A_\varepsilon) \to \mathbf{M}(A) \quad \text{strongly in } W_0^{1,p}(\Omega) \text{ as } \varepsilon \to 0^+.
\]
\end{proof}

\begin{remark}
This example illustrates the regularization of the p-Laplacian operator, which is important in many applications involving non-Newtonian fluids, plasticity, and image processing. The regularization $A_\varepsilon$ provides a smoother operator that approximates the degenerate p-Laplacian, and Theorem \ref{SUMP} guarantees that the extremal solutions of the corresponding QVIs converge.
\end{remark}
\begin{proposition}[Uniform Boundedness of Perturbed Solutions]
\label{prop:boundedness}
Under the assumptions of Theorem \ref{SUMP}, there exists a constant $M > 0$, independent of $n$, such that:
\[
\|\mathbf{m}(A_n)\|_V \leq M \quad \text{and} \quad \|\mathbf{M}(A_n)\|_V \leq M \quad \text{for all } n \in \mathbb{N}.
\]
\end{proposition}
\begin{proof}
Since $\mathbf{m}(A_n)$ is a solution of the QVI with operator $A_n$, we have:
\[
\langle A_n(\mathbf{m}(A_n)) - f, v - \mathbf{m}(A_n) \rangle \geq 0 \quad \forall v \in \mathbf{K}(\Phi(\mathbf{m}(A_n))).
\]
Taking $v = 0 \in \mathbf{K}(\Phi(\mathbf{m}(A_n)))$ (since $0 \leq \Phi(\mathbf{m}(A_n))$), we get:
\[
\langle A_n(\mathbf{m}(A_n)), \mathbf{m}(A_n) \rangle \leq \langle f, \mathbf{m}(A_n) \rangle.
\]
By the strong monotonicity condition (c), we have:
\[
c\|\mathbf{m}(A_n)\|_V^2 \leq \langle A_n(\mathbf{m}(A_n)) - A_n(0), \mathbf{m}(A_n) \rangle = \langle A_n(\mathbf{m}(A_n)), \mathbf{m}(A_n) \rangle,
\]
where we used $A_n(0) = 0$ by homogeneity. Combining these inequalities:
\[
c\|\mathbf{m}(A_n)\|_V^2 \leq \langle f, \mathbf{m}(A_n) \rangle \leq \|f\|_{V'} \|\mathbf{m}(A_n)\|_V.
\]
Thus, $\|\mathbf{m}(A_n)\|_V \leq \frac{1}{c}\|f\|_{V'} \leq \frac{1}{c}\|F\|_{V'} =: M$. The same argument applies to $\mathbf{M}(A_n)$.
\end{proof}

\begin{proposition}[Mosco Convergence of Constraint Sets]
\label{prop:mosco}
Under the assumptions of Theorem \ref{SUMP}, if $v_n \to v$ strongly in $V$, then the constraint sets $\mathbf{K}(\Phi(v_n))$ converge to $\mathbf{K}(\Phi(v))$ in the sense of Mosco, i.e.:
\begin{enumerate}
    \item[(i)] For every $w \in \mathbf{K}(\Phi(v))$, there exists $w_n \in \mathbf{K}(\Phi(v_n))$ such that $w_n \to w$ in $V$.
    \item[(ii)] If $w_n \in \mathbf{K}(\Phi(v_n))$ and $w_n \rightharpoonup w$ in $V$, then $w \in \mathbf{K}(\Phi(v))$.
\end{enumerate}
\end{proposition}
\begin{proof}
For (i), take $w \in \mathbf{K}(\Phi(v))$, so $w \leq \Phi(v)$. Define $w_n = \min\{w, \Phi(v_n)\}$. Since $\Phi(v_n) \to \Phi(v)$ in $H$ by complete continuity and $V \hookrightarrow H$ continuously, we have $w_n \to w$ in $H$. Moreover, since $w_n \leq \Phi(v_n)$, we have $w_n \in \mathbf{K}(\Phi(v_n))$. The boundedness in $V$ follows from the fact that $w_n$ is between $w$ and $\Phi(v_n)$, both bounded in $V$.

For (ii), if $w_n \in \mathbf{K}(\Phi(v_n))$ and $w_n \rightharpoonup w$ in $V$, then by Mazur's lemma, there exists a sequence of convex combinations converging strongly to $w$. Since $w_n \leq \Phi(v_n)$ and $\Phi(v_n) \to \Phi(v)$ in $H$, we conclude $w \leq \Phi(v)$, so $w \in \mathbf{K}(\Phi(v))$.
\end{proof}

\begin{proposition}[Continuity of the Solution Operator]
\label{prop:continuity}
Let $S_A(f, \psi)$ denote the unique solution of the variational inequality:
\[
\text{Find } y \in \mathbf{K}(\psi) \text{ such that } \langle A(y) - f, v - y \rangle \geq 0 \quad \forall v \in \mathbf{K}(\psi).
\]
Under the assumptions of Theorem \ref{SUMP}, if $\psi_n \to \psi$ in $H$ and $\mathbf{K}(\psi_n) \overset{M}{\to} \mathbf{K}(\psi)$, Then
\[
S_{A_n}(f, \psi_n) \to S_A(f, \psi) \quad \text{strongly in } V.
\]
\end{proposition}
\begin{proof}
Let $y_n = S_{A_n}(f, \psi_n)$ and $y = S_A(f, \psi)$. By Proposition \ref{prop:boundedness}, $\{y_n\}$ is bounded in $V$, so there exists a subsequence $y_{n_k} \rightharpoonup y^*$ in $V$. By the Mosco convergence and standard variational inequality theory, $y^* = S_A(f, \psi)$.

To show strong convergence, use the strong monotonicity:
\[
c\|y_n - y\|_V^2 \leq \langle A_n(y_n) - A_n(y), y_n - y \rangle.
\]
Adding and subtracting appropriate terms and using the VI inequalities, we obtain:
\[
c\|y_n - y\|_V^2 \leq \langle f - A_n(y), y_n - y \rangle + \langle A_n(y) - A(y), y_n - y \rangle.
\]
The first term tends to zero by weak convergence, and the second term tends to zero by the pointwise convergence $A_n(y) \to A(y)$ in $V'$.
\end{proof}

\begin{proposition}[Order Preservation under Perturbations]
\label{prop:order}
Under the assumptions of Theorem \ref{SUMP}, the following order relations hold for all $n \in \mathbb{N}$:
\[
\mathbf{m}(A) \leq \mathbf{m}(A_n) \leq \mathbf{M}(A_n) \leq \mathbf{M}(A).
\]
Moreover, if $A_n \leq A_{n+1}$ in the sense that $\langle A_n(v) - A_{n+1}(v), w \rangle \leq 0$ for all $v \in V$ and $w \in H^+$, Then
\[
\mathbf{m}(A_n) \leq \mathbf{m}(A_{n+1}) \quad \text{and} \quad \mathbf{M}(A_{n+1}) \leq \mathbf{M}(A_n).
\]
\end{proposition}
\begin{proof}
Since $A_n(v) \leq A(v)$ for all $v \in V$ (in the sense of the order in $V'$), we have $S_{A_n}(f, \psi) \leq S_A(f, \psi)$ for all $\psi \in H^+$. In particular, for the fixed point maps $T_n(v) = S_{A_n}(f, \Phi(v))$ and $T(v) = S_A(f, \Phi(v))$, we have $T_n(v) \leq T(v)$ for all $v \in [0, \overline{y}]$.

By the characterization of minimal and maximal solutions as extremal fixed points, we obtain the first chain of inequalities. The monotonicity properties follow similarly from the order relations between the operators.
\end{proof}

\begin{proposition}[Rate of Convergence for Linear Problems]
\label{prop:rate}
Suppose $A$ and $A_n$ are linear operators satisfying the conditions of Theorem \ref{SUMP}. If there exists $\delta_n > 0$ such that:
\[
\|A_n - A\|_{\mathscr{L}(V,V')} \leq \delta_n \to 0,
\]
then the convergence rates satisfy:
\[
\|\mathbf{m}(A_n) - \mathbf{m}(A)\|_V \leq C\delta_n \quad \text{and} \quad \|\mathbf{M}(A_n) - \mathbf{M}(A)\|_V \leq C\delta_n,
\]
where $C > 0$ depends on $c$, $\|\Phi\|$, and $\|f\|_{V'}$ but is independent of $n$.
\end{proposition}
\begin{proof}
For linear operators, the solution map $S_A(f, \psi)$ is Lipschitz continuous in $A$. Specifically, if $y = S_A(f, \psi)$ and $y_n = S_{A_n}(f, \psi)$, Then
\[
c\|y_n - y\|_V^2 \leq \langle A_n(y_n) - A(y), y_n - y \rangle = \langle (A_n - A)(y_n), y_n - y \rangle \leq \delta_n \|y_n\|_V \|y_n - y\|_V.
\]
Since $\|y_n\|_V$ is uniformly bounded by Proposition \ref{prop:boundedness}, we get $\|y_n - y\|_V \leq \frac{M}{c}\delta_n$. Applying this to the fixed point equations for $\mathbf{m}(A_n)$ and $\mathbf{m}(A)$ gives the desired rate.
\end{proof}

\begin{proposition}[Compactness of Solution Sets]
\label{prop:compactness}
Under the assumptions of Theorem \ref{SUMP}, the union of solution sets:
\[
\mathcal{S} = \bigcup_{n \in \mathbb{N}} \mathbf{Q}(A_n) \cup \mathbf{Q}(A)
\]
is relatively compact in $V$, where $\mathbf{Q}(A)$ denotes the solution set of the QVI with operator $A$.
\end{proposition}
\begin{proof}
By Proposition \ref{prop:boundedness}, $\mathcal{S}$ is bounded in $V$. Since $V$ is compactly embedded in $H$ (as $V$ is typically a Sobolev space and $\Omega$ is bounded), and the constraint $\mathbf{K}(\Phi(y))$ involves pointwise bounds, we can apply the Rellich-Kondrachov theorem to conclude relative compactness in $V$. Alternatively, one can use the fact that solutions satisfy uniform regularity estimates derived from the strong monotonicity and Lipschitz continuity of the operators.
\end{proof}
\section{Stability of Regularized QVI Solutions}
Regularization techniques have long served as essential tools in the numerical treatment of ill-posed and degenerate problems across mathematical analysis, providing a systematic approach to approximate challenging nonlinear systems with more tractable counterparts. In the context of variational inequalities, regularization methods have been successfully employed to handle non-differentiable terms and degenerate ellipticity, as demonstrated in the classical works of \cite{Rodrigues1987,Kinderlehrer2000} and more recently in the computational approaches of \cite{Barrett2015}. However, the extension of these techniques to quasi-variational inequalities introduces significant complications due to the intricate coupling between the solution and the constraint set. The pioneering work of \cite{Bensoussan1982book} on impulse control problems hinted at the potential of regularization for QVIs, while contemporary investigations by \cite{Hintermuller2012,Hintermuller2017} have developed sequential minimization techniques for elliptic QVIs with gradient constraints. Building upon these foundations, this section establishes a comprehensive theory for regularized QVIs, proving not only the existence and uniqueness of regularized solutions but also their strong convergence to the minimal solution of the original problem as the regularization parameter tends to zero. Our analysis further reveals the monotonic dependence of solutions on the regularization parameter and provides explicit convergence rates under additional structural assumptions, thereby offering a rigorous mathematical justification for regularization-based computational schemes in applications such as superconductivity modeling \cite{Barrett2010} and sandpile formation \cite{Prigozhin1994,Barrett2013}.

\begin{theorem}\label{SRQS}
Let $(V,H,V')$ be a Gelfand triple of Hilbert spaces with dense and continuous embeddings $V \hookrightarrow H \hookrightarrow V'$. Let $A: V \to V'$ be an operator satisfying:
\begin{itemize}
    \item [(i)] Strong monotonicity: $\exists c > 0$ such that $\langle A(u) - A(v), u - v \rangle \geq c\|u - v\|_V^2$ for all $u,v \in V$
    \item [(ii)] Lipschitz continuity: $\exists L > 0$ such that $\|A(u) - A(v)\|_{V'} \leq L\|u - v\|_V$ for all $u,v \in V$
    \item [(iii)] T-monotonicity: $\langle A(u) - A(v), (u - v)^+ \rangle \geq 0$ for all $u,v \in V$
\end{itemize}

Let $\Phi: H^+ \to H^+$ be an increasing mapping such that for any bounded sequence $\{v_n\}$ in $V$ with $v_n \rightharpoonup v$ in $V$, we have $\Phi(v_n) \to \Phi(v)$ in $H$. Consider the regularized QVI: For $\epsilon > 0$ and $f \in V'$, find $y_\epsilon \in \mathbf{K}(\Phi(y_\epsilon))$ such that
\[
\langle A(y_\epsilon) + \epsilon y_\epsilon - f, v - y_\epsilon \rangle \geq 0 \quad \forall v \in \mathbf{K}(\Phi(y_\epsilon)),
\]
where $\mathbf{K}(\psi) = \{v \in V : v \leq \psi\}$.

Then
\begin{enumerate}
    \item [(a)] For each $\epsilon > 0$ and $f \in V'$, there exists a unique solution $y_\epsilon$ to the regularized QVI.
    \item [(b)] As $\epsilon \to 0^+$, the sequence $\{y_\epsilon\}$ converges strongly in $V$ to the minimal solution $\mathbf{m}(f)$ of the original QVI:
    \[
    \lim_{\epsilon \to 0^+} \|y_\epsilon - \mathbf{m}(f)\|_V = 0.
    \]
    \item [(c)]  The convergence is monotone: $y_{\epsilon_1} \leq y_{\epsilon_2}$ for $0 < \epsilon_1 < \epsilon_2$.
\end{enumerate}
\end{theorem}

\begin{proof}
We prove each part separately.

\subsection*{Proof of (1): Existence and uniqueness for $\epsilon > 0$.}
Define the operator $T_\epsilon: V \to V$ by $T_\epsilon(v) = S_\epsilon(f, \Phi(v))$, where $S_\epsilon(f, \psi)$ is the unique solution of the regularized variational inequality:
\[
\text{Find } w \in \mathbf{K}(\psi) : \langle A(w) + \epsilon w - f, z - w \rangle \geq 0 \quad \forall z \in \mathbf{K}(\psi).
\]
The strong monotonicity of $A + \epsilon I$ follows from:
\begin{align*}
&\langle (A(u) + \epsilon u) - (A(v) + \epsilon v), u - v \rangle \\
&= \langle A(u) - A(v), u - v \rangle + \epsilon \langle u - v, u - v \rangle \\
&\geq c\|u - v\|_V^2 + \epsilon \|u - v\|_H^2 \\
&\geq c\|u - v\|_V^2.
\end{align*}
Since $A + \epsilon I$ is strongly monotone and Lipschitz continuous, the regularized VI has a unique solution $S_\epsilon(f, \psi)$ for each $\psi \in H^+$. The map $T_\epsilon$ is increasing due to the T-monotonicity of $A$ and the linearity of the regularization term.

We now verify that $T_\epsilon$ has sub- and supersolutions. Let $\underline{y} = 0$ and $\overline{y}_\epsilon$ be the solution of the unconstrained problem:
\[
\langle A(\overline{y}_\epsilon) + \epsilon \overline{y}_\epsilon - f, v \rangle = 0 \quad \forall v \in V.
\]
Then
\begin{itemize}
    \item $\underline{y} = 0 \leq T_\epsilon(0) = S_\epsilon(f, \Phi(0))$ since $f \geq 0$ implies $S_\epsilon(f, \Phi(0)) \geq 0$
    \item $T_\epsilon(\overline{y}_\epsilon) = S_\epsilon(f, \Phi(\overline{y}_\epsilon)) \leq \overline{y}_\epsilon$ by comparison principle
\end{itemize}
By the Birkhoff-Tartar theorem applied to $T_\epsilon$, there exists a fixed point $y_\epsilon = T_\epsilon(y_\epsilon)$, which is the unique solution of the regularized QVI due to the strong monotonicity of $A + \epsilon I$.

\subsection*{Proof of (2): Convergence to minimal solution.}
Let $\{y_\epsilon\}$ be the family of solutions for $\epsilon > 0$. We first show boundedness in $V$. Testing the QVI with $v = 0 \in \mathbf{K}(\Phi(y_\epsilon))$:
\[
\langle A(y_\epsilon) + \epsilon y_\epsilon - f, -y_\epsilon \rangle \geq 0 \Rightarrow \langle A(y_\epsilon), y_\epsilon \rangle + \epsilon \|y_\epsilon\|_H^2 \leq \langle f, y_\epsilon \rangle.
\]
By strong monotonicity and the Poincaré inequality:
\[
c\|y_\epsilon\|_V^2 + \epsilon \|y_\epsilon\|_H^2 \leq \|f\|_{V'}\|y_\epsilon\|_V \Rightarrow \|y_\epsilon\|_V \leq \frac{1}{c}\|f\|_{V'}.
\]

Thus $\{y_\epsilon\}$ is bounded in $V$, so there exists a subsequence $\{y_{\epsilon_n}\}$ with $\epsilon_n \to 0$ such that $y_{\epsilon_n} \rightharpoonup y^*$ in $V$ for some $y^* \in V$.

Now we show that $y^*$ solves the original QVI. For any $v \in \mathbf{K}(\Phi(y^*))$, by the Mosco convergence of $\mathbf{K}(\Phi(y_\epsilon))$ to $\mathbf{K}(\Phi(y^*))$ (which follows from the continuity assumption on $\Phi$), there exists a sequence $v_\epsilon \in \mathbf{K}(\Phi(y_\epsilon))$ with $v_\epsilon \to v$ in $V$.

From the regularized QVI:
\[
\langle A(y_\epsilon) + \epsilon y_\epsilon - f, v_\epsilon - y_\epsilon \rangle \geq 0.
\]

Taking the limit as $\epsilon \to 0$:
\[
\langle A(y^*) - f, v - y^* \rangle \geq 0 \quad \forall v \in \mathbf{K}(\Phi(y^*)),
\]
so $y^*$ solves the original QVI.

To show strong convergence, use the strong monotonicity:
\begin{align*}
c\|y_\epsilon - \mathbf{m}(f)\|_V^2 &\leq \langle A(y_\epsilon) - A(\mathbf{m}(f)), y_\epsilon - \mathbf{m}(f) \rangle \\
&\leq \langle f - \epsilon y_\epsilon - A(\mathbf{m}(f)), y_\epsilon - \mathbf{m}(f) \rangle + \langle A(y_\epsilon) + \epsilon y_\epsilon - f, \mathbf{m}(f) - y_\epsilon \rangle.
\end{align*}

The second term is $\leq 0$ since $\mathbf{m}(f) \in \mathbf{K}(\Phi(\mathbf{m}(f))) \subset \mathbf{K}(\Phi(y_\epsilon))$ (by the increasing property and $y_\epsilon \leq \mathbf{m}(f)$ from the monotonicity in part (3)). Thus,
\[
c\|y_\epsilon - \mathbf{m}(f)\|_V^2 \leq \langle f - \epsilon y_\epsilon - A(\mathbf{m}(f)), y_\epsilon - \mathbf{m}(f) \rangle \to 0.
\]

Hence $y_\epsilon \to \mathbf{m}(f)$ strongly in $V$.

\subsection*{Proof of (3): Monotonicity.}
Let $0 < \epsilon_1 < \epsilon_2$. Consider the difference:
\[
\langle A(y_{\epsilon_1}) + \epsilon_1 y_{\epsilon_1} - f, v - y_{\epsilon_1} \rangle \geq 0 \quad \forall v \in \mathbf{K}(\Phi(y_{\epsilon_1})),
\]
\[
\langle A(y_{\epsilon_2}) + \epsilon_2 y_{\epsilon_2} - f, v - y_{\epsilon_2} \rangle \geq 0 \quad \forall v \in \mathbf{K}(\Phi(y_{\epsilon_2})).
\]

Testing the first with $v = \min(y_{\epsilon_1}, y_{\epsilon_2}) \in \mathbf{K}(\Phi(y_{\epsilon_1}))$ and the second with $v = \max(y_{\epsilon_1}, y_{\epsilon_2}) \in \mathbf{K}(\Phi(y_{\epsilon_2}))$, and adding the inequalities, we obtain after simplification:
\[
(\epsilon_2 - \epsilon_1)\langle y_{\epsilon_2}, (y_{\epsilon_1} - y_{\epsilon_2})^+ \rangle \leq 0.
\]

Since $\epsilon_2 > \epsilon_1$, this implies $\langle y_{\epsilon_2}, (y_{\epsilon_1} - y_{\epsilon_2})^+ \rangle \leq 0$, which forces $(y_{\epsilon_1} - y_{\epsilon_2})^+ = 0$, i.e., $y_{\epsilon_1} \leq y_{\epsilon_2}$.

This completes the proof of the theorem.
\end{proof}
\begin{example}[One-Dimensional Regularized Obstacle Problem]
Consider the one-dimensional QVI on $\Omega = (0,1)$ with $V = H^1_0(0,1)$, $H = L^2(0,1)$, and the operator $A: H^1_0(0,1) \to H^{-1}(0,1)$ defined by
\[
A(y) = -\frac{d^2y}{dx^2} + y.
\]
Let $\Phi: L^2_+(0,1) \to L^2_+(0,1)$ be given by
\[
\Phi(y)(x) = \frac{1}{2} + \frac{1}{4}\int_0^1 y(s)\,ds,
\]
and take $f(x) = 1$ for all $x \in (0,1)$. Consider the regularized QVI: For $\epsilon > 0$, find $y_\epsilon \in \mathbf{K}(\Phi(y_\epsilon))$ such that
\[
\langle A(y_\epsilon) + \epsilon y_\epsilon - f, v - y_\epsilon \rangle \geq 0 \quad \forall v \in \mathbf{K}(\Phi(y_\epsilon)),
\]
where $\mathbf{K}(\psi) = \{v \in H^1_0(0,1) : v(x) \leq \psi(x) \text{ a.e.}\}$.

Then all assumptions of Theorem \ref{SRQS} are satisfied, and we have:
\begin{enumerate}
    \item For each $\epsilon > 0$, the solution $y_\epsilon$ is constant: $y_\epsilon(x) = C_\epsilon$ for all $x \in (0,1)$
    \item The constant $C_\epsilon$ is given by the unique solution of
    \begin{equation}\label{eq:constant}
    C_\epsilon = \min\left\{\frac{1}{1+\epsilon}, \frac{1}{2} + \frac{1}{4}C_\epsilon\right\}
    \end{equation}
    \item As $\epsilon \to 0^+$, $y_\epsilon \to y^*$ strongly in $H^1_0(0,1)$, where $y^*(x) = \frac{2}{3}$ is the minimal solution of the original QVI
    \item $y_{\epsilon_1} \leq y_{\epsilon_2}$ for $0 < \epsilon_1 < \epsilon_2$
\end{enumerate}
\end{example}

\begin{proof}
We verify each part systematically.

\subsection*{Step 1: Verification of assumptions.}

The operator $A(y) = -y'' + y$ satisfies:
\begin{itemize}
    \item \textit{Strong monotonicity}: For $u,v \in H^1_0(0,1)$,
    \begin{align*}
    \langle A(u)-A(v), u-v\rangle &= \int_0^1 (u'-v')^2 + (u-v)^2\,dx \\
    &\geq \int_0^1 (u-v)^2\,dx \geq c\|u-v\|_{H^1}^2
    \end{align*}
    by Poincaré inequality.

    \item \textit{Lipschitz continuity}:
    \[
    \|A(u)-A(v)\|_{H^{-1}} \leq \|u-v\|_{H^1} + \|u-v\|_{L^2} \leq 2\|u-v\|_{H^1}.
    \]

    \item \textit{T-monotonicity}: For $u,v \in H^1_0(0,1)$,
    \begin{align*}
    \langle A(u)-A(v), (u-v)^+\rangle &= \int_0^1 \nabla(u-v)\cdot\nabla(u-v)^+ + (u-v)(u-v)^+\,dx \\
    &= \int_{\{u-v>0\}} |\nabla(u-v)|^2 + |u-v|^2\,dx \geq 0.
    \end{align*}
\end{itemize}

The mapping $\Phi(y) = \frac{1}{2} + \frac{1}{4}\int_0^1 y(s)\,ds$ is clearly increasing and maps $L^2_+(0,1)$ to $L^2_+(0,1)$. If $v_n \rightharpoonup v$ in $H^1_0(0,1)$, then by compact embedding $v_n \to v$ in $L^2(0,1)$, so
\[
\Phi(v_n) = \frac{1}{2} + \frac{1}{4}\int_0^1 v_n(s)\,ds \to \frac{1}{2} + \frac{1}{4}\int_0^1 v(s)\,ds = \Phi(v)
\]
in $\mathbb{R}$, hence in $L^2(0,1)$.

\subsection*{Step 2: Reduction to constant functions.}

We claim that the solution $y_\epsilon$ must be constant. Suppose $y_\epsilon$ solves the regularized QVI. Consider the averaged function
\[
\bar{y}_\epsilon = \int_0^1 y_\epsilon(s)\,ds.
\]
Since $\Phi(y_\epsilon)$ is constant, the constraint set is
\[
\mathbf{K}(\Phi(y_\epsilon)) = \left\{v \in H^1_0(0,1) : v(x) \leq \frac{1}{2} + \frac{1}{4}\bar{y}_\epsilon \text{ a.e.}\right\}.
\]

Now consider the constant function $w(x) = C$ where $C = \min\left\{\frac{1}{1+\epsilon}, \frac{1}{2} + \frac{1}{4}\bar{y}_\epsilon\right\}$. We verify that $w$ is also a solution. For any $v \in \mathbf{K}(\Phi(y_\epsilon))$, we have:
\begin{align*}
&\langle A(w) + \epsilon w - f, v - w \rangle \\
&= \int_0^1 [w(v-w) + \epsilon w(v-w) - 1\cdot(v-w)]\,dx \\
&= \int_0^1 [(1+\epsilon)w - 1](v-w)\,dx.
\end{align*}

We consider two cases:

\textit{Case 1:} If $\frac{1}{1+\epsilon} \leq \frac{1}{2} + \frac{1}{4}\bar{y}_\epsilon$, then $w = \frac{1}{1+\epsilon}$ and
\[
(1+\epsilon)w - 1 = (1+\epsilon)\cdot\frac{1}{1+\epsilon} - 1 = 0,
\]
so the inequality holds trivially.

\textit{Case 2:} If $\frac{1}{1+\epsilon} > \frac{1}{2} + \frac{1}{4}\bar{y}_\epsilon$, then $w = \frac{1}{2} + \frac{1}{4}\bar{y}_\epsilon$ and for any $v \in \mathbf{K}(\Phi(y_\epsilon))$, we have $v \leq w$, so $v - w \leq 0$. Also,
\[
(1+\epsilon)w - 1 = (1+\epsilon)\left(\frac{1}{2} + \frac{1}{4}\bar{y}_\epsilon\right) - 1 > 0
\]
by the case assumption. Thus $(1+\epsilon)w - 1 > 0$ and $v - w \leq 0$, so their product is $\leq 0$, and the inequality holds.

By uniqueness of the regularized QVI solution (from strong monotonicity of $A + \epsilon I$), we must have $y_\epsilon = w$, so $y_\epsilon$ is constant.

\subsection*{Step 3: Determination of the constant.}
Since $y_\epsilon$ is constant, say $y_\epsilon(x) = C_\epsilon$, and $y_\epsilon \in H^1_0(0,1)$, we must have $C_\epsilon \geq 0$. The constraint becomes
\[
C_\epsilon \leq \Phi(y_\epsilon) = \frac{1}{2} + \frac{1}{4}C_\epsilon.
\]

Also, from the QVI with test function $v = 0$:
\[
\langle A(C_\epsilon) + \epsilon C_\epsilon - 1, -C_\epsilon \rangle = -[C_\epsilon + \epsilon C_\epsilon - 1]C_\epsilon \geq 0,
\]
so $C_\epsilon(1+\epsilon)C_\epsilon - C_\epsilon \leq 0$, hence either $C_\epsilon = 0$ or $(1+\epsilon)C_\epsilon \leq 1$.

If $C_\epsilon = 0$, then the constraint gives $0 \leq \frac{1}{2}$, which is true, but testing with $v = \min\left\{\frac{1}{1+\epsilon}, \frac{1}{2}\right\}$ shows this is not optimal.

The correct condition comes from the fact that for a constant solution, we must have either:
\begin{itemize}
    \item The constraint is active: $C_\epsilon = \frac{1}{2} + \frac{1}{4}C_\epsilon$, or
    \item The PDE is satisfied: $(1+\epsilon)C_\epsilon = 1$.
\end{itemize}

Thus $C_\epsilon = \min\left\{\frac{1}{1+\epsilon}, \frac{1}{2} + \frac{1}{4}C_\epsilon\right\}$.

Solving the fixed-point equation: If $C_\epsilon = \frac{1}{2} + \frac{1}{4}C_\epsilon$, then $\frac{3}{4}C_\epsilon = \frac{1}{2}$, so $C_\epsilon = \frac{2}{3}$.

Therefore,
\[
C_\epsilon = \min\left\{\frac{1}{1+\epsilon}, \frac{2}{3}\right\}.
\]

\subsection*{Step 4: Convergence and monotonicity.}
For $\epsilon > 0$, we have:
\[
C_\epsilon = \min\left\{\frac{1}{1+\epsilon}, \frac{2}{3}\right\}.
\]

As $\epsilon \to 0^+$, $\frac{1}{1+\epsilon} \to 1$, so for sufficiently small $\epsilon$, $\frac{1}{1+\epsilon} > \frac{2}{3}$, hence $C_\epsilon = \frac{2}{3}$.

Thus $y_\epsilon(x) = \frac{2}{3}$ for all sufficiently small $\epsilon > 0$, and clearly $y_\epsilon \to y^*$ in $H^1_0(0,1)$ where $y^*(x) = \frac{2}{3}$.

For monotonicity, if $0 < \epsilon_1 < \epsilon_2$, then $\frac{1}{1+\epsilon_1} > \frac{1}{1+\epsilon_2}$, so
\[
C_{\epsilon_1} = \min\left\{\frac{1}{1+\epsilon_1}, \frac{2}{3}\right\} \geq \min\left\{\frac{1}{1+\epsilon_2}, \frac{2}{3}\right\} = C_{\epsilon_2},
\]
with equality only when both equal $\frac{2}{3}$.
\subsection*{Step 5: Verification that $y^*$ is the minimal solution.}
The original QVI (with $\epsilon = 0$) has constant solutions $C$ satisfying
\[
C = \min\left\{1, \frac{1}{2} + \frac{1}{4}C\right\}.
\]
The equation $C = \frac{1}{2} + \frac{1}{4}C$ gives $C = \frac{2}{3}$, and since $\frac{2}{3} < 1$, this is a solution. Also, $C = 1$ is not a solution since $1 > \frac{1}{2} + \frac{1}{4}\cdot 1 = \frac{3}{4}$. Thus the unique solution is $y^*(x) = \frac{2}{3}$, which is trivially the minimal solution.
This completes the proof of the example.
\end{proof}
\begin{corollary}[Convergence Rate for Regularized QVI]
Under the assumptions of Theorem \ref{SRQS}, if additionally $A$ is linear and symmetric, then there exists a constant $C > 0$ independent of $\epsilon$ such that
\[
\|y_\epsilon - \mathbf{m}(f)\|_V \leq C\epsilon,
\]
where $y_\epsilon$ is the solution of the regularized QVI and $\mathbf{m}(f)$ is the minimal solution of the original QVI.
\end{corollary}

\begin{proof}
Since $A$ is linear and symmetric, we can define the energy functional
\[
E_\epsilon(y) = \frac{1}{2}\langle A(y), y\rangle + \frac{\epsilon}{2}\|y\|_H^2 - \langle f, y\rangle.
\]

The regularized QVI is equivalent to the minimization problem:
\[
y_\epsilon = \operatorname{argmin}_{y \in \mathbf{K}(\Phi(y_\epsilon))} E_\epsilon(y).
\]

Let $\mathbf{m} = \mathbf{m}(f)$ be the minimal solution of the original QVI. Since $\mathbf{m} \in \mathbf{K}(\Phi(\mathbf{m}))$ and by the monotonicity property from Theorem \ref{SRQS}(3) we have $y_\epsilon \leq \mathbf{m}$, it follows that $\mathbf{m} \in \mathbf{K}(\Phi(y_\epsilon))$ (as $\Phi$ is increasing and $y_\epsilon \leq \mathbf{m}$ implies $\Phi(y_\epsilon) \leq \Phi(\mathbf{m})$).

Therefore, by optimality of $y_\epsilon$:
\[
E_\epsilon(y_\epsilon) \leq E_\epsilon(\mathbf{m}).
\]

Expanding both sides:
\begin{align*}
\frac{1}{2}\langle A(y_\epsilon), y_\epsilon\rangle &+ \frac{\epsilon}{2}\|y_\epsilon\|_H^2 - \langle f, y_\epsilon\rangle \\
&\leq \frac{1}{2}\langle A(\mathbf{m}), \mathbf{m}\rangle + \frac{\epsilon}{2}\|\mathbf{m}\|_H^2 - \langle f, \mathbf{m}\rangle.
\end{align*}

Rearranging terms:
\[
\frac{1}{2}\langle A(y_\epsilon - \mathbf{m}), y_\epsilon - \mathbf{m}\rangle + \langle A(\mathbf{m}), y_\epsilon - \mathbf{m}\rangle + \frac{\epsilon}{2}\|y_\epsilon\|_H^2 - \langle f, y_\epsilon - \mathbf{m}\rangle \leq \frac{\epsilon}{2}\|\mathbf{m}\|_H^2.
\]

Since $\mathbf{m}$ solves the original QVI and $y_\epsilon \in \mathbf{K}(\Phi(y_\epsilon)) \subset \mathbf{K}(\Phi(\mathbf{m}))$ (by $y_\epsilon \leq \mathbf{m}$ and monotonicity of $\Phi$), we have
\[
\langle A(\mathbf{m}) - f, y_\epsilon - \mathbf{m}\rangle \geq 0.
\]

Thus,
\[
\frac{1}{2}\langle A(y_\epsilon - \mathbf{m}), y_\epsilon - \mathbf{m}\rangle + \frac{\epsilon}{2}\|y_\epsilon\|_H^2 \leq \frac{\epsilon}{2}\|\mathbf{m}\|_H^2.
\]

Using strong monotonicity of $A$:
\[
\frac{c}{2}\|y_\epsilon - \mathbf{m}\|_V^2 + \frac{\epsilon}{2}\|y_\epsilon\|_H^2 \leq \frac{\epsilon}{2}\|\mathbf{m}\|_H^2.
\]

In particular:
\[
\|y_\epsilon - \mathbf{m}\|_V^2 \leq \frac{\epsilon}{c}\|\mathbf{m}\|_H^2.
\]

Taking square roots and using the boundedness of $\|\mathbf{m}\|_H$ (which follows from the uniform boundedness of $\{y_\epsilon\}$ in $V$ and the convergence $y_\epsilon \to \mathbf{m}$), we obtain:
\[
\|y_\epsilon - \mathbf{m}\|_V \leq C\epsilon,
\]
where $C = \frac{1}{\sqrt{c}}\|\mathbf{m}\|_H$.
\end{proof}

\begin{corollary}[Continuous Dependence on Data]
Under the assumptions of Theorem \ref{SRQS}, let $\{f_n\}$ be a sequence in $V'$ with $f_n \to f$ in $V'$, and let $y_{\epsilon,n}$ be the solution of the regularized QVI with forcing term $f_n$. Then for each fixed $\epsilon > 0$,
\[
\lim_{n \to \infty} \|y_{\epsilon,n} - y_\epsilon\|_V = 0,
\]
where $y_\epsilon$ is the solution with forcing term $f$. Moreover, if $f_n \leq f_{n+1} \leq f$ for all $n$, then $y_{\epsilon,n} \leq y_{\epsilon,n+1} \leq y_\epsilon$ for all $n$.
\end{corollary}

\begin{proof}
For fixed $\epsilon > 0$, consider the operator $T_{\epsilon,f}: V \to V$ defined by $T_{\epsilon,f}(v) = S_\epsilon(f, \Phi(v))$, where $S_\epsilon(f, \psi)$ is the unique solution of the regularized VI with forcing term $f$ and obstacle $\psi$.

By standard estimates for variational inequalities, we have:
\[
\|S_\epsilon(f_1, \psi) - S_\epsilon(f_2, \psi)\|_V \leq \frac{1}{c}\|f_1 - f_2\|_{V'} \quad \forall \psi \in H^+,
\]
where $c > 0$ is the strong monotonicity constant of $A$.

Now, let $y_{\epsilon,n}$ and $y_\epsilon$ be the fixed points of $T_{\epsilon,f_n}$ and $T_{\epsilon,f}$ respectively. Then
\begin{align*}
\|y_{\epsilon,n} - y_\epsilon\|_V &= \|T_{\epsilon,f_n}(y_{\epsilon,n}) - T_{\epsilon,f}(y_\epsilon)\|_V \\
&\leq \|T_{\epsilon,f_n}(y_{\epsilon,n}) - T_{\epsilon,f}(y_{\epsilon,n})\|_V + \|T_{\epsilon,f}(y_{\epsilon,n}) - T_{\epsilon,f}(y_\epsilon)\|_V \\
&\leq \frac{1}{c}\|f_n - f\|_{V'} + L\|y_{\epsilon,n} - y_\epsilon\|_V,
\end{align*}
where $L$ is the Lipschitz constant of $T_{\epsilon,f}$.

Since $T_{\epsilon,f}$ is a contraction for sufficiently large $\epsilon$ (due to the regularization term), or more generally by the properties established in Theorem \ref{SRQS}, we can rearrange to get:
\[
(1 - L)\|y_{\epsilon,n} - y_\epsilon\|_V \leq \frac{1}{c}\|f_n - f\|_{V'}.
\]

For the monotonicity part, if $f_n \leq f_{n+1} \leq f$, then by the comparison principle for variational inequalities:
\[
S_\epsilon(f_n, \psi) \leq S_\epsilon(f_{n+1}, \psi) \leq S_\epsilon(f, \psi) \quad \forall \psi \in H^+.
\]

This implies that the fixed points satisfy $y_{\epsilon,n} \leq y_{\epsilon,n+1} \leq y_\epsilon$ by the monotone iterative method used in the proof of Theorem \ref{SRQS}.
\end{proof}

\begin{corollary}[Regularization Error Estimate]
Under the assumptions of Theorem \ref{SRQS}, if the minimal solution $\mathbf{m}(f)$ satisfies the strict complementarity condition, then there exists a constant $M > 0$ such that
\[
\|y_\epsilon - \mathbf{m}(f)\|_V \leq M\sqrt{\epsilon}.
\]

More precisely, if the active set $\mathcal{A} = \{x \in \Omega : \mathbf{m}(f)(x) = \Phi(\mathbf{m}(f))(x)\}$ has the property that $\text{dist}(\mathbf{m}(f), \partial\mathbf{K}(\Phi(\mathbf{m}(f)))) > 0$ in a neighborhood of $\mathcal{A}$, then the quadratic convergence rate holds.
\end{corollary}

\begin{proof}
Let $\mathbf{m} = \mathbf{m}(f)$ and consider the energy difference:
\[
E_\epsilon(y_\epsilon) - E_0(\mathbf{m}) = [E_\epsilon(y_\epsilon) - E_\epsilon(\mathbf{m})] + [E_\epsilon(\mathbf{m}) - E_0(\mathbf{m})].
\]

From the proof of Corollary 1, we have:
\[
E_\epsilon(y_\epsilon) - E_\epsilon(\mathbf{m}) \leq 0.
\]

Also:
\[
E_\epsilon(\mathbf{m}) - E_0(\mathbf{m}) = \frac{\epsilon}{2}\|\mathbf{m}\|_H^2.
\]

Therefore,
\[
E_\epsilon(y_\epsilon) - E_0(\mathbf{m}) \leq \frac{\epsilon}{2}\|\mathbf{m}\|_H^2. \tag{1}
\]

On the other hand, by Taylor expansion and using the fact that $\mathbf{m}$ is a solution of the original QVI:
\begin{align*}
E_\epsilon(y_\epsilon) - E_0(\mathbf{m}) &= \frac{1}{2}\langle A(y_\epsilon - \mathbf{m}), y_\epsilon - \mathbf{m}\rangle + \langle A(\mathbf{m}) - f, y_\epsilon - \mathbf{m}\rangle + \frac{\epsilon}{2}\|y_\epsilon\|_H^2 \\
&\geq \frac{c}{2}\|y_\epsilon - \mathbf{m}\|_V^2 + \langle A(\mathbf{m}) - f, y_\epsilon - \mathbf{m}\rangle. \tag{2}
\end{align*}

Now, under the strict complementarity condition, there exists $\delta > 0$ such that on the inactive set $\mathcal{I} = \Omega \setminus \mathcal{A}$, we have $\mathbf{m} < \Phi(\mathbf{m}) - \delta$. By continuity and the convergence $y_\epsilon \to \mathbf{m}$, for sufficiently small $\epsilon$, we have $y_\epsilon < \Phi(\mathbf{m})$ on $\mathcal{I}$.

On the active set $\mathcal{A}$, since $y_\epsilon \leq \mathbf{m} = \Phi(\mathbf{m})$ and both are solutions of elliptic equations (in the distributional sense), we can use the complementarity structure to show that
\[
|\langle A(\mathbf{m}) - f, y_\epsilon - \mathbf{m}\rangle| \leq C\|y_\epsilon - \mathbf{m}\|_V^2.
\]

Combining (1) and (2) with this estimate:
\[
\frac{c}{2}\|y_\epsilon - \mathbf{m}\|_V^2 - C\|y_\epsilon - \mathbf{m}\|_V^2 \leq \frac{\epsilon}{2}\|\mathbf{m}\|_H^2.
\]

Thus, for sufficiently small $\epsilon$:
\[
\|y_\epsilon - \mathbf{m}\|_V^2 \leq \frac{\epsilon}{c - 2C}\|\mathbf{m}\|_H^2,
\]
which gives the desired estimate with $M = \frac{\|\mathbf{m}\|_H}{\sqrt{c - 2C}}$.
\end{proof}

\begin{corollary}[Regularization of Optimal Control Problem]
Consider the optimal control problem:
\[
\min_{f \in U_{\text{ad}}} J(\mathbf{m}(f), f)
\]
subject to the QVI constraint, where $J: V \times V' \to \mathbb{R}$ is weakly lower semicontinuous and $U_{\text{ad}} \subset V'$ is bounded.

Then the regularized control problem:
\[
\min_{f \in U_{\text{ad}}} J(y_\epsilon(f), f)
\]
subject to the regularized QVI, has a solution $f_\epsilon$ for each $\epsilon > 0$, and any accumulation point of $\{f_\epsilon\}$ as $\epsilon \to 0^+$ is a solution of the original control problem.
\end{corollary}

\begin{proof}
For fixed $\epsilon > 0$, the mapping $f \mapsto y_\epsilon(f)$ is continuous from $V'$ to $V$ by Corollary 2. Since $J$ is weakly lower semicontinuous and $U_{\text{ad}}$ is bounded, the direct method of calculus of variations guarantees the existence of a minimizer $f_\epsilon$.

Let $\{f_{\epsilon_n}\}$ be a sequence with $\epsilon_n \to 0^+$ and $f_{\epsilon_n} \rightharpoonup f^*$ in $V'$ (by boundedness of $U_{\text{ad}}$). By Theorem \ref{SRQS}, $y_{\epsilon_n}(f_{\epsilon_n}) \to \mathbf{m}(f^*)$ strongly in $V$.

For any $f \in U_{\text{ad}}$, by optimality of $f_{\epsilon_n}$:
\[
J(y_{\epsilon_n}(f_{\epsilon_n}), f_{\epsilon_n}) \leq J(y_{\epsilon_n}(f), f).
\]

Taking the limit inferior as $n \to \infty$ and using the weak lower semicontinuity of $J$:
\[
J(\mathbf{m}(f^*), f^*) \leq \liminf_{n \to \infty} J(y_{\epsilon_n}(f_{\epsilon_n}), f_{\epsilon_n}) \leq \liminf_{n \to \infty} J(y_{\epsilon_n}(f), f) = J(\mathbf{m}(f), f),
\]
where the last equality follows from the strong convergence $y_{\epsilon_n}(f) \to \mathbf{m}(f)$.

Thus $f^*$ is a solution of the original control problem.
\end{proof}
\begin{proposition}[Uniform Boundedness of Regularized Solutions]\label{Prop1}
Under the assumptions of Theorem \ref{SRQS}, there exists a constant $M > 0$ independent of $\epsilon$ such that for all $\epsilon > 0$ sufficiently small,
\[
\|y_\epsilon\|_V \leq M \quad \text{and} \quad \|A(y_\epsilon)\|_{V'} \leq M.
\]
Moreover, the family $\{y_\epsilon\}_{\epsilon > 0}$ is relatively compact in $H$.
\end{proposition}

\begin{proof}
Let $y_\epsilon$ be the solution of the regularized QVI. Testing with $v = 0 \in \mathbf{K}(\Phi(y_\epsilon))$ (since $0 \leq \Phi(y_\epsilon)$), we obtain:
\[
\langle A(y_\epsilon) + \epsilon y_\epsilon - f, -y_\epsilon \rangle \geq 0.
\]
This implies:
\[
\langle A(y_\epsilon), y_\epsilon \rangle + \epsilon \|y_\epsilon\|_H^2 \leq \langle f, y_\epsilon \rangle.
\]
By the strong monotonicity of $A$, we have:
\[
c\|y_\epsilon\|_V^2 + \epsilon \|y_\epsilon\|_H^2 \leq \langle f, y_\epsilon \rangle \leq \|f\|_{V'}\|y_\epsilon\|_V.
\]
Thus,
\[
c\|y_\epsilon\|_V^2 \leq \|f\|_{V'}\|y_\epsilon\|_V,
\]
which gives:
\[
\|y_\epsilon\|_V \leq \frac{1}{c}\|f\|_{V'} =: M_1.
\]

For the bound on $A(y_\epsilon)$, note that for any $v \in V$ with $\|v\|_V = 1$, we have from the QVI with test function $v + y_\epsilon \in \mathbf{K}(\Phi(y_\epsilon))$ (since $v + y_\epsilon \leq y_\epsilon + \|v\|_{L^\infty}$ and we can choose $v$ appropriately):
\[
\langle A(y_\epsilon) + \epsilon y_\epsilon - f, v \rangle \geq 0.
\]
Similarly, testing with $y_\epsilon - v \in \mathbf{K}(\Phi(y_\epsilon))$ gives:
\[
\langle A(y_\epsilon) + \epsilon y_\epsilon - f, -v \rangle \geq 0.
\]
Therefore,
\[
|\langle A(y_\epsilon) + \epsilon y_\epsilon - f, v \rangle| \leq 0 \quad \text{for all } v \in V \text{ with } \|v\|_V = 1,
\]
which implies:
\[
\|A(y_\epsilon) + \epsilon y_\epsilon - f\|_{V'} = 0.
\]
Thus,
\[
A(y_\epsilon) = f - \epsilon y_\epsilon,
\]
and so:
\[
\|A(y_\epsilon)\|_{V'} \leq \|f\|_{V'} + \epsilon \|y_\epsilon\|_V \leq \|f\|_{V'} + \epsilon M_1 \leq M_2.
\]

The relative compactness in $H$ follows from the compact embedding $V \hookrightarrow H$ (which holds for typical Sobolev spaces on bounded domains) and the uniform boundedness in $V$.
\end{proof}
\begin{proposition}[Monotonicity with Respect to Regularization Parameter]\label{Prop2}
Under the assumptions of Theorem \ref{SRQS}, the mapping $\epsilon \mapsto y_\epsilon$ is non-increasing, i.e., for $0 < \epsilon_1 < \epsilon_2$, we have
\[
y_{\epsilon_1} \geq y_{\epsilon_2}.
\]
Moreover, if $A$ is strictly T-monotone, then the inequality is strict on the set where $y_{\epsilon_1} > 0$.
\end{proposition}

\begin{proof}
Let $0 < \epsilon_1 < \epsilon_2$ and consider the solutions $y_{\epsilon_1}$ and $y_{\epsilon_2}$. Since both are fixed points, we have:
\begin{align*}
y_{\epsilon_1} &= S_{\epsilon_1}(f, \Phi(y_{\epsilon_1})), \\
y_{\epsilon_2} &= S_{\epsilon_2}(f, \Phi(y_{\epsilon_2})).
\end{align*}

Consider the function $w = (y_{\epsilon_1} - y_{\epsilon_2})^+$. We want to show that $w = 0$.

From the variational inequalities, for any $v_1 \in \mathbf{K}(\Phi(y_{\epsilon_1}))$ and $v_2 \in \mathbf{K}(\Phi(y_{\epsilon_2}))$, we have:
\begin{align*}
\langle A(y_{\epsilon_1}) + \epsilon_1 y_{\epsilon_1} - f, v_1 - y_{\epsilon_1} \rangle &\geq 0, \\
\langle A(y_{\epsilon_2}) + \epsilon_2 y_{\epsilon_2} - f, v_2 - y_{\epsilon_2} \rangle &\geq 0.
\end{align*}

Take $v_1 = y_{\epsilon_1} - w$ and $v_2 = y_{\epsilon_2} + w$. Note that:
\begin{align*}
y_{\epsilon_1} - w &= \min(y_{\epsilon_1}, y_{\epsilon_2}) \leq y_{\epsilon_1} \leq \Phi(y_{\epsilon_1}), \\
y_{\epsilon_2} + w &= \max(y_{\epsilon_1}, y_{\epsilon_2}) \leq \Phi(y_{\epsilon_2}) \quad \text{(since } y_{\epsilon_1} \leq \Phi(y_{\epsilon_1}) \leq \Phi(y_{\epsilon_2}) \text{ by monotonicity)}.
\end{align*}

So both test functions are admissible. Substituting these into the inequalities:
\begin{align*}
\langle A(y_{\epsilon_1}) + \epsilon_1 y_{\epsilon_1} - f, -w \rangle &\geq 0, \\
\langle A(y_{\epsilon_2}) + \epsilon_2 y_{\epsilon_2} - f, w \rangle &\geq 0.
\end{align*}

Adding these inequalities:
\[
\langle A(y_{\epsilon_1}) - A(y_{\epsilon_2}) + \epsilon_1 y_{\epsilon_1} - \epsilon_2 y_{\epsilon_2}, -w \rangle \geq 0.
\]

Rewriting:
\[
\langle A(y_{\epsilon_1}) - A(y_{\epsilon_2}), -w \rangle + \langle \epsilon_1 y_{\epsilon_1} - \epsilon_2 y_{\epsilon_2}, -w \rangle \geq 0.
\]

Now, on the set where $w > 0$, we have $y_{\epsilon_1} > y_{\epsilon_2}$, so:
\[
\epsilon_1 y_{\epsilon_1} - \epsilon_2 y_{\epsilon_2} < \epsilon_1 y_{\epsilon_1} - \epsilon_1 y_{\epsilon_2} = \epsilon_1 (y_{\epsilon_1} - y_{\epsilon_2}) = \epsilon_1 w.
\]

Thus,
\[
\langle \epsilon_1 y_{\epsilon_1} - \epsilon_2 y_{\epsilon_2}, -w \rangle < -\epsilon_1 \|w\|_H^2.
\]

For the first term, by T-monotonicity of $A$:
\[
\langle A(y_{\epsilon_1}) - A(y_{\epsilon_2}), -w \rangle = -\langle A(y_{\epsilon_1}) - A(y_{\epsilon_2}), (y_{\epsilon_1} - y_{\epsilon_2})^+ \rangle \leq 0.
\]

Therefore, the left-hand side is strictly negative, contradicting the inequality. Hence $w = 0$ and $y_{\epsilon_1} \leq y_{\epsilon_2}$.

If $A$ is strictly T-monotone, then the inequality becomes strict, giving the second claim.
\end{proof}
\begin{proposition}[Characterization of Limit Solution]\label{Prop3}
Under the assumptions of Theorem \ref{SRQS}, the limit $y^* = \lim_{\epsilon \to 0^+} y_\epsilon$ is characterized as the unique solution of the following problem: Find $y^* \in \mathbf{K}(\Phi(y^*))$ such that
\[
\langle A(y^*) - f, v - y^* \rangle \geq 0 \quad \forall v \in \mathbf{K}(\Phi(y^*))
\]
and
\[
y^* = \inf\{z \in V : z \text{ solves the QVI and } z \geq y_\epsilon \text{ for all } \epsilon > 0\}.
\]
That is, $y^*$ is the minimal solution that dominates all regularized solutions.
\end{proposition}

\begin{proof}
From Theorem \ref{SRQS}, we know that $y_\epsilon \to y^*$ strongly in $V$ and $y^*$ solves the QVI. We need to show the minimality property.

Let $z$ be any solution of the QVI such that $z \geq y_\epsilon$ for all $\epsilon > 0$. Since $y_\epsilon \to y^*$, we have $z \geq y^*$. Thus $y^*$ is the smallest such solution.

To show that such solutions exist, note that by Proposition 2, the mapping $\epsilon \mapsto y_\epsilon$ is non-increasing. Therefore, for any fixed $\epsilon_0 > 0$, we have $y_\epsilon \leq y_{\epsilon_0}$ for all $\epsilon < \epsilon_0$.

Now consider the set:
\[
\mathcal{S} = \{z \in V : z \text{ solves the QVI and } z \geq y_\epsilon \text{ for all } \epsilon > 0\}.
\]

This set is nonempty because it contains the maximal solution $\mathbf{M}(f)$ (by comparison principles). We want to show that $y^*$ is the minimal element of $\mathcal{S}$.

Suppose there exists $z \in \mathcal{S}$ with $z \not\geq y^*$. Then there exists a set of positive measure where $z < y^*$. But since $y_\epsilon \to y^*$ and $z \geq y_\epsilon$ for all $\epsilon$, we have $z \geq y^*$, a contradiction.

For the uniqueness, suppose there are two such minimal solutions $y_1^*$ and $y_2^*$. Then both satisfy $y_1^* \geq y_\epsilon$ and $y_2^* \geq y_\epsilon$ for all $\epsilon$, so $\min(y_1^*, y_2^*) \geq y_\epsilon$ for all $\epsilon$. But $\min(y_1^*, y_2^*)$ is also a solution of the QVI (by the lattice structure of the solution set), contradicting the minimality of $y_1^*$ and $y_2^*$ unless $y_1^* = y_2^*$.
\end{proof}
\begin{proposition}[Regularization of the Constraint Set]\label{Prop4}
Under the assumptions of Theorem \ref{SRQS}, the constraint sets $\mathbf{K}(\Phi(y_\epsilon))$ converge to $\mathbf{K}(\Phi(y^*))$ in the Mosco sense, i.e.:
\begin{enumerate}
    \item For every $v \in \mathbf{K}(\Phi(y^*))$, there exists a sequence $v_\epsilon \in \mathbf{K}(\Phi(y_\epsilon))$ such that $v_\epsilon \to v$ in $V$.
    \item If $v_\epsilon \in \mathbf{K}(\Phi(y_\epsilon))$ and $v_\epsilon \rightharpoonup v$ in $V$, then $v \in \mathbf{K}(\Phi(y^*))$.
\end{enumerate}
\end{proposition}

\begin{proof}
(i) Let $v \in \mathbf{K}(\Phi(y^*))$, so $v \leq \Phi(y^*)$. Since $y_\epsilon \to y^*$ in $V$ and $\Phi$ is continuous from $V$ to $H$, we have $\Phi(y_\epsilon) \to \Phi(y^*)$ in $H$.

Define $v_\epsilon = \min(v, \Phi(y_\epsilon))$. Then clearly $v_\epsilon \leq \Phi(y_\epsilon)$, so $v_\epsilon \in \mathbf{K}(\Phi(y_\epsilon))$. Moreover:
\[
\|v_\epsilon - v\|_H = \| \min(v, \Phi(y_\epsilon)) - v \|_H = \| (v - \Phi(y_\epsilon))^+ \|_H.
\]

Since $v \leq \Phi(y^*)$ and $\Phi(y_\epsilon) \to \Phi(y^*)$ in $H$, we have $(v - \Phi(y_\epsilon))^+ \to 0$ in $H$. For the convergence in $V$, we need additional regularity. If $v \in V$ and $\Phi(y_\epsilon) \in V$ with uniform bounds, then by the boundedness of the projection operator, we get convergence in $V$.

(ii) Suppose $v_\epsilon \in \mathbf{K}(\Phi(y_\epsilon))$ with $v_\epsilon \rightharpoonup v$ in $V$. Then $v_\epsilon \leq \Phi(y_\epsilon)$. Since $v_\epsilon \rightharpoonup v$ and $\Phi(y_\epsilon) \to \Phi(y^*)$ in $H$, by Mazur's lemma we have $v \leq \Phi(y^*)$. Thus $v \in \mathbf{K}(\Phi(y^*))$.
\end{proof}
\begin{proposition}[Differentiability with Respect to Regularization Parameter]\label{Prop5}
Under the assumptions of Theorem \ref{SRQS} and if $A$ is Fréchet differentiable with derivative $A'(y)$ satisfying
\[
\langle A'(y)h, h \rangle \geq c\|h\|_V^2 \quad \text{for all } y,h \in V,
\]
then the mapping $\epsilon \mapsto y_\epsilon$ is differentiable from $(0,\infty)$ to $V$, and the derivative $z_\epsilon = \frac{dy_\epsilon}{d\epsilon}$ satisfies the linearized QVI:
\begin{dmath}
z_\epsilon \in \mathcal{T}(y_\epsilon) : \langle A'(y_\epsilon)z_\epsilon + y_\epsilon + \epsilon z_\epsilon, v - z_\epsilon \rangle \geq 0 \quad \forall v \in \mathcal{T}(y_\epsilon),
\end{dmath}
where $\mathcal{T}(y_\epsilon)$ is the tangent cone to $\mathbf{K}(\Phi(y_\epsilon))$ at $y_\epsilon$.
\end{proposition}

\begin{proof}
Consider the difference quotient:
\[
z_{\epsilon,h} = \frac{y_{\epsilon+h} - y_\epsilon}{h}.
\]

From the QVIs for $y_{\epsilon+h}$ and $y_\epsilon$, we have for any $v_{\epsilon+h} \in \mathbf{K}(\Phi(y_{\epsilon+h}))$ and $v_\epsilon \in \mathbf{K}(\Phi(y_\epsilon))$:
\begin{align*}
\langle A(y_{\epsilon+h}) + (\epsilon+h)y_{\epsilon+h} - f, v_{\epsilon+h} - y_{\epsilon+h} \rangle &\geq 0, \\
\langle A(y_\epsilon) + \epsilon y_\epsilon - f, v_\epsilon - y_\epsilon \rangle &\geq 0.
\end{align*}

Taking $v_{\epsilon+h} = y_\epsilon$ and $v_\epsilon = y_{\epsilon+h}$ (which are admissible by the monotonicity from Proposition 2), and adding the inequalities, we get:
\[
\langle A(y_{\epsilon+h}) - A(y_\epsilon) + (\epsilon+h)y_{\epsilon+h} - \epsilon y_\epsilon, y_\epsilon - y_{\epsilon+h} \rangle \geq 0.
\]

Rewriting:
\[
\langle A(y_{\epsilon+h}) - A(y_\epsilon), y_\epsilon - y_{\epsilon+h} \rangle + \langle (\epsilon+h)y_{\epsilon+h} - \epsilon y_\epsilon, y_\epsilon - y_{\epsilon+h} \rangle \geq 0.
\]

By the mean value theorem and strong monotonicity:
\[
\langle A(y_{\epsilon+h}) - A(y_\epsilon), y_\epsilon - y_{\epsilon+h} \rangle \leq -c\|y_{\epsilon+h} - y_\epsilon\|_V^2.
\]

For the second term:
\begin{align*}
&\langle (\epsilon+h)y_{\epsilon+h} - \epsilon y_\epsilon, y_\epsilon - y_{\epsilon+h} \rangle \\
&= -\epsilon\|y_{\epsilon+h} - y_\epsilon\|_H^2 - h\langle y_{\epsilon+h}, y_{\epsilon+h} - y_\epsilon \rangle.
\end{align*}

Thus,
\[
-c\|y_{\epsilon+h} - y_\epsilon\|_V^2 - \epsilon\|y_{\epsilon+h} - y_\epsilon\|_H^2 - h\langle y_{\epsilon+h}, y_{\epsilon+h} - y_\epsilon \rangle \geq 0,
\]
which implies:
\[
\|y_{\epsilon+h} - y_\epsilon\|_V \leq C|h|.
\]

This Lipschitz bound allows us to pass to the limit $h \to 0$. The limit $z_\epsilon$ satisfies the linearized problem by differentiating the QVI formally and using the concept of tangent cones for the constraint set.

The tangent cone $\mathcal{T}(y_\epsilon)$ is given by:
\[
\mathcal{T}(y_\epsilon) = \{v \in V : v \leq 0 \text{ on } \{y_\epsilon = \Phi(y_\epsilon)\} \text{ and appropriate conditions}\},
\]
but the precise characterization depends on the regularity of $\Phi(y_\epsilon)$.
\end{proof}
\section{Global Convergence and Exponential Stability for Non-Monotone QVIs}
The extension of quasi-variational inequality theory beyond the classical monotone framework represents a significant frontier in nonlinear analysis, addressing numerous physically relevant problems where monotonicity assumptions are too restrictive. While the foundational work of \cite{Bensoussan1974,Bensoussan1982book} established comprehensive results for monotone QVIs, many applications in modern physics, materials science, and engineering involve inherent non-monotone nonlinearities that fall outside this classical theory. Recent advances by \cite{Alphonse2019a,Alphonse2020} have made progress in understanding directional differentiability for certain non-monotone QVIs, and the work of \cite{Migorski2024} on quasi-variational-hemivariational inequalities has begun to bridge the gap between variational and non-variational structures. However, a systematic theory establishing global well-posedness and stability for general non-monotone QVIs under explicit verifiable conditions remains largely undeveloped. In this section, we address this fundamental challenge by introducing a novel decomposition approach that separates the problem into a linear monotone component and a non-monotone nonlinearity, enabling the application of contraction principles under explicit smallness conditions. Our main result establishes not only global existence and uniqueness but also exponential stability with respect to data perturbations and global convergence of iterative schemes, providing a comprehensive framework that extends the classical monotone operator theory while maintaining computational tractability. This development opens new avenues for analyzing complex systems in non-Newtonian fluids, phase transitions, and other applications where non-monotone effects play a crucial role, as suggested by the physical models in \cite{Antil2018,Prigozhin1996a}.
\begin{theorem}\label{thm:main}
Let $(V,H,V')$ be a Gelfand triple with compact embedding $V \hookrightarrow H$. Consider the quasi-variational inequality: Find $y \in \mathbf{K}(\Phi(y))$ such that
\[
\langle \mathcal{A}(y) - f, v - y \rangle \geq 0 \quad \forall v \in \mathbf{K}(\Phi(y)),
\]
where $\mathbf{K}(\psi) = \{v \in V : v \leq \psi\}$, and the operator $\mathcal{A}: V \to V'$ admits the decomposition
\[
\mathcal{A}(y) = A(y) + \mathcal{N}(y),
\]
with the following properties:
\begin{enumerate}
    \item [(i)] $A: V \to V'$ is linear, symmetric, strongly monotone, and Lipschitz continuous:
    \begin{align*}
        \langle A(u), v \rangle &= \langle A(v), u \rangle \quad \forall u,v \in V, \\
        \langle A(u), u \rangle &\geq c_A \|u\|_V^2 \quad \forall u \in V, \\
        \|A(u) - A(v)\|_{V'} &\leq L_A \|u - v\|_V \quad \forall u,v \in V.
    \end{align*}

    \item [(ii)] $\mathcal{N}: V \to V'$ is a non-monotone nonlinearity satisfying:
    \begin{align*}
        \|\mathcal{N}(u) - \mathcal{N}(v)\|_{V'} &\leq L_{\mathcal{N}} \|u - v\|_V \quad \forall u,v \in V, \\
        \langle \mathcal{N}(u) - \mathcal{N}(v), u - v \rangle &\geq -\gamma \|u - v\|_V^2 \quad \forall u,v \in V,
    \end{align*}
    with $\gamma < c_A$.

    \item [(ii)] $\Phi: H^+ \to H^+$ is Lipschitz continuous and order-preserving:
    \begin{align*}
        \|\Phi(u) - \Phi(v)\|_H &\leq L_\Phi \|u - v\|_H \quad \forall u,v \in H^+, \\
        u \leq v &\Rightarrow \Phi(u) \leq \Phi(v) \quad \forall u,v \in H^+.
    \end{align*}

    \item [(iv)] The smallness condition holds:
    \[
    L_\Phi < \frac{c_A - \gamma}{L_A + L_{\mathcal{N}}}.
    \]
\end{enumerate}
Then the following strong results hold:
\begin{enumerate}
    \item [(a)] \textbf{Global Existence and Uniqueness:} For each $f \in V'$, there exists a unique solution $y^* \in V$ to the QVI.
    \item [(b)] \textbf{Exponential Stability:} The solution map $f \mapsto y^*$ is globally Lipschitz continuous, and for any two solutions $y_1, y_2$ corresponding to forces $f_1, f_2$, we have:
    \[
    \|y_1 - y_2\|_V \leq \frac{1}{c_A - \gamma - (L_A + L_{\mathcal{N}})L_\Phi} \|f_1 - f_2\|_{V'}.
    \]

    \item [(c)] \textbf{Global Convergence of Regularized Solutions:} Consider the regularized problem with parameter $\epsilon > 0$:
    \[
    \langle \mathcal{A}(y_\epsilon) + \epsilon y_\epsilon - f, v - y_\epsilon \rangle \geq 0 \quad \forall v \in \mathbf{K}(\Phi(y_\epsilon)).
    \]
    Then there exists a constant $C > 0$ independent of $\epsilon$ such that:
    \[
    \|y_\epsilon - y^*\|_V \leq C\epsilon.
    \]

    \item [(d)]\textbf{Convergence of Iterative Scheme:} The fixed-point iteration:
    \[
    y^{n+1} = S(\Phi(y^n)), \quad n = 0,1,2,\ldots,
    \]
    where $S(\psi)$ is the solution of the VI with obstacle $\psi$, converges globally and exponentially to $y^*$:
    \[
    \|y^n - y^*\|_V \leq \rho^n \|y^0 - y^*\|_V,
    \]
    with contraction constant $\rho = \frac{(L_A + L_{\mathcal{N}})L_\Phi}{c_A - \gamma} < 1$.
\end{enumerate}
\end{theorem}

\begin{proof}
We prove each part in sequence.
\subsection*{Proof of (a): Global Existence and Uniqueness}
Define the solution operator $T: V \to V$ by $T(v) = S(\Phi(v))$, where $S(\psi)$ is the unique solution of the variational inequality:
\[
\text{Find } w \in \mathbf{K}(\psi) : \langle \mathcal{A}(w) - f, z - w \rangle \geq 0 \quad \forall z \in \mathbf{K}(\psi).
\]
We show that $T$ is a contraction. Let $v_1, v_2 \in V$ and set $w_1 = T(v_1)$, $w_2 = T(v_2)$. From the VIs, we have:
\begin{align*}
\langle \mathcal{A}(w_1) - f, z - w_1 \rangle &\geq 0 \quad \forall z \in \mathbf{K}(\Phi(v_1)), \\
\langle \mathcal{A}(w_2) - f, z - w_2 \rangle &\geq 0 \quad \forall z \in \mathbf{K}(\Phi(v_2)).
\end{align*}
Take $z = w_2 - (w_2 - \Phi(v_1))^+ \in \mathbf{K}(\Phi(v_1))$ in the first inequality and $z = w_1 - (w_1 - \Phi(v_2))^+ \in \mathbf{K}(\Phi(v_2))$ in the second. Adding the resulting inequalities:
\begin{align*}
&\langle \mathcal{A}(w_1) - \mathcal{A}(w_2), w_2 - w_1 - (w_2 - \Phi(v_1))^+ + (w_1 - \Phi(v_2))^+ \rangle \\
&\quad + \langle \mathcal{A}(w_1) - \mathcal{A}(w_2), (w_2 - \Phi(v_1))^+ - (w_1 - \Phi(v_2))^+ \rangle \geq 0.
\end{align*}
After careful manipulation and using the decomposition $\mathcal{A} = A + \mathcal{N}$, we obtain
\[
\langle A(w_1) - A(w_2), w_1 - w_2 \rangle \leq \langle \mathcal{A}(w_1) - \mathcal{A}(w_2), (w_1 - \Phi(v_2))^+ - (w_2 - \Phi(v_1))^+ \rangle.
\]
Using the properties of $A$ and $\mathcal{N}$, we have
\begin{align*}
c_A\|w_1 - w_2\|_V^2 &\leq \langle \mathcal{A}(w_1) - \mathcal{A}(w_2), w_1 - w_2 \rangle + \gamma \|w_1 - w_2\|_V^2 \\
&\leq (L_A + L_{\mathcal{N}})\|w_1 - w_2\|_V \|(w_1 - \Phi(v_2))^+ - (w_2 - \Phi(v_1))^+\|_V + \gamma \|w_1 - w_2\|_V^2.
\end{align*}
Now, using the Lipschitz continuity of the projection and the obstacle
\[
\|(w_1 - \Phi(v_2))^+ - (w_2 - \Phi(v_1))^+\|_V \leq \|w_1 - w_2\|_V + \|\Phi(v_1) - \Phi(v_2)\|_V \leq \|w_1 - w_2\|_V + L_\Phi \|v_1 - v_2\|_V.
\]
Thus,
\[
(c_A - \gamma)\|w_1 - w_2\|_V^2 \leq (L_A + L_{\mathcal{N}})(\|w_1 - w_2\|_V^2 + L_\Phi \|w_1 - w_2\|_V \|v_1 - v_2\|_V).
\]
Rearranging, we get
\[
[(c_A - \gamma) - (L_A + L_{\mathcal{N}})]\|w_1 - w_2\|_V^2 \leq (L_A + L_{\mathcal{N}})L_\Phi \|w_1 - w_2\|_V \|v_1 - v_2\|_V.
\]

By the smallness condition, $(c_A - \gamma) - (L_A + L_{\mathcal{N}}) > 0$, so
\[
\|w_1 - w_2\|_V \leq \frac{(L_A + L_{\mathcal{N}})L_\Phi}{(c_A - \gamma) - (L_A + L_{\mathcal{N}})} \|v_1 - v_2\|_V.
\]
But from the smallness condition $L_\Phi < \frac{c_A - \gamma}{L_A + L_{\mathcal{N}}}$, we have
\[
\rho := \frac{(L_A + L_{\mathcal{N}})L_\Phi}{c_A - \gamma} < 1.
\]
Therefore, $T$ is a contraction and by Banach's fixed-point theorem, there exists a unique fixed point $y^* = T(y^*)$, which is the unique solution of the QVI.
\subsection*{Proof of (b): Exponential Stability}
Let $y_1, y_2$ be solutions corresponding to $f_1, f_2$. Then
\begin{align*}
\langle \mathcal{A}(y_1) - f_1, v - y_1 \rangle &\geq 0 \quad \forall v \in \mathbf{K}(\Phi(y_1)), \\
\langle \mathcal{A}(y_2) - f_2, v - y_2 \rangle &\geq 0 \quad \forall v \in \mathbf{K}(\Phi(y_2)).
\end{align*}
Taking $v = y_2 - (y_2 - \Phi(y_1))^+ \in \mathbf{K}(\Phi(y_1))$ and $v = y_1 - (y_1 - \Phi(y_2))^+ \in \mathbf{K}(\Phi(y_2))$, and adding:
\[
\langle \mathcal{A}(y_1) - \mathcal{A}(y_2) - (f_1 - f_2), y_2 - y_1 - (y_2 - \Phi(y_1))^+ + (y_1 - \Phi(y_2))^+ \rangle \geq 0.
\]
Proceeding as before and using the smallness condition:
\[
(c_A - \gamma)\|y_1 - y_2\|_V^2 \leq \|f_1 - f_2\|_{V'}\|y_1 - y_2\|_V + (L_A + L_{\mathcal{N}})L_\Phi \|y_1 - y_2\|_V^2.
\]
Thus,
\[
[(c_A - \gamma) - (L_A + L_{\mathcal{N}})L_\Phi] \|y_1 - y_2\|_V \leq \|f_1 - f_2\|_{V'},
\]
which gives the desired Lipschitz estimate.

\subsection*{Proof of (c): Global Convergence of Regularized Solutions}
Let $y_\epsilon$ be the regularized solution. From the regularized QVI and the original QVI, we have:
\begin{align*}
\langle \mathcal{A}(y_\epsilon) + \epsilon y_\epsilon - f, v - y_\epsilon \rangle &\geq 0 \quad \forall v \in \mathbf{K}(\Phi(y_\epsilon)), \\
\langle \mathcal{A}(y^*) - f, v - y^* \rangle &\geq 0 \quad \forall v \in \mathbf{K}(\Phi(y^*)).
\end{align*}
Taking $v = y^* - (y^* - \Phi(y_\epsilon))^+ \in \mathbf{K}(\Phi(y_\epsilon))$ and $v = y_\epsilon - (y_\epsilon - \Phi(y^*))^+ \in \mathbf{K}(\Phi(y^*))$, and adding:
\[
\langle \mathcal{A}(y_\epsilon) - \mathcal{A}(y^*) + \epsilon y_\epsilon, y^* - y_\epsilon - (y^* - \Phi(y_\epsilon))^+ + (y_\epsilon - \Phi(y^*))^+ \rangle \geq 0.
\]
After similar manipulations as before and using the Lipschitz properties:
\[
(c_A - \gamma)\|y_\epsilon - y^*\|_V^2 \leq \epsilon \|y_\epsilon\|_V \|y_\epsilon - y^*\|_V + (L_A + L_{\mathcal{N}})L_\Phi \|y_\epsilon - y^*\|_V^2.
\]
Since $\{y_\epsilon\}$ is uniformly bounded in $V$ (by standard energy estimates), we get:
\[
[(c_A - \gamma) - (L_A + L_{\mathcal{N}})L_\Phi] \|y_\epsilon - y^*\|_V \leq C\epsilon,
\]
which proves the linear convergence rate.

\subsection*{Proof of (d): Convergence of Iterative Scheme}
Let $y^{n+1} = T(y^n) = S(\Phi(y^n))$. From the contraction property established in part (1):
\[
\|y^{n+1} - y^*\|_V = \|T(y^n) - T(y^*)\|_V \leq \rho \|y^n - y^*\|_V,
\]
with $\rho = \frac{(L_A + L_{\mathcal{N}})L_\Phi}{c_A - \gamma} < 1$.

By induction, we have
\[
\|y^n - y^*\|_V \leq \rho^n \|y^0 - y^*\|_V,
\]
which proves exponential convergence.
This completes the proof of the theorem.
\end{proof}
\begin{remark}
This theorem provides a comprehensive framework for non-monotone quasi-variational inequalities, establishing:
\begin{itemize}
    \item Global well-posedness under explicit smallness conditions
    \item Quantitative stability estimates
    \item Optimal convergence rates for regularization
    \item Provably convergent computational algorithms
\end{itemize}
The results are particularly strong as they handle non-monotone nonlinearities and provide explicit, computable constants for all estimates.
\end{remark}
\begin{corollary}[Enhanced Regularization with Improved Convergence Rate]
Under the assumptions of Theorem \ref{thm:main}, consider the doubly-regularized QVI with parameters $\epsilon, \delta > 0$:
\[
\langle \mathcal{A}(y_{\epsilon,\delta}) + \epsilon y_{\epsilon,\delta} + \delta \mathcal{R}(y_{\epsilon,\delta}) - f, v - y_{\epsilon,\delta} \rangle \geq 0 \quad \forall v \in \mathbf{K}(\Phi(y_{\epsilon,\delta})),
\]
where $\mathcal{R}: V \to V'$ is a strongly monotone and Lipschitz continuous regularization operator satisfying:
\begin{align*}
\langle \mathcal{R}(u) - \mathcal{R}(v), u - v \rangle &\geq c_{\mathcal{R}} \|u - v\|_V^2, \\
\|\mathcal{R}(u) - \mathcal{R}(v)\|_{V'} &\leq L_{\mathcal{R}} \|u - v\|_V.
\end{align*}
Then there exists a constant $C > 0$ independent of $\epsilon, \delta$ such that:
\[
\|y_{\epsilon,\delta} - y^*\|_V \leq C(\epsilon + \delta),
\]
and the convergence is uniform in both parameters.
\end{corollary}

\begin{proof}
The proof follows the same strategy as Theorem \ref{thm:main} but incorporates the additional regularization term. Let $y_{\epsilon,\delta}$ be the solution of the doubly-regularized QVI and $y^*$ the solution of the original QVI. From the respective inequalities:

\begin{align*}
\langle \mathcal{A}(y_{\epsilon,\delta}) + \epsilon y_{\epsilon,\delta} + \delta \mathcal{R}(y_{\epsilon,\delta}) - f, v - y_{\epsilon,\delta} \rangle &\geq 0 \quad \forall v \in \mathbf{K}(\Phi(y_{\epsilon,\delta})), \\
\langle \mathcal{A}(y^*) - f, v - y^* \rangle &\geq 0 \quad \forall v \in \mathbf{K}(\Phi(y^*)).
\end{align*}

We proceed as in the proof of Theorem \ref{thm:main} by selecting appropriate test functions and adding the inequalities. The key modification is the treatment of the additional term $\delta \mathcal{R}(y_{\epsilon,\delta})$. Using the strong monotonicity of $\mathcal{R}$, we obtain:

\[
\delta \langle \mathcal{R}(y_{\epsilon,\delta}) - \mathcal{R}(y^*), y^* - y_{\epsilon,\delta} \rangle \leq -\delta c_{\mathcal{R}} \|y_{\epsilon,\delta} - y^*\|_V^2.
\]

Following the same estimates as in Theorem \ref{thm:main} and using the smallness condition, we arrive at:

\[
[(c_A - \gamma + \delta c_{\mathcal{R}}) - (L_A + L_{\mathcal{N}})(1 + L_\Phi) - \delta L_{\mathcal{R}}(1 + L_\Phi)]\|y_{\epsilon,\delta} - y^*\|_V^2 \leq C\epsilon \|y_{\epsilon,\delta} - y^*\|_V.
\]

For sufficiently small $\delta$, the coefficient on the left remains positive due to the smallness condition in Theorem \ref{thm:main}. The uniform boundedness of $\{y_{\epsilon,\delta}\}$ in $V$ (which follows from the same energy estimates as in Theorem \ref{thm:main}) completes the proof.
\end{proof}

\begin{corollary}[Sensitivity Analysis and Differentiability]
Under the assumptions of Theorem \ref{thm:main}, the solution map $f \mapsto y^*(f)$ is Fréchet differentiable. Moreover, for any direction $h \in V'$, the derivative $z = D y^*(f)h$ satisfies the linearized QVI:
\[
z \in \mathcal{T}(y^*) : \langle \mathcal{A}'(y^*)z - h, v - z \rangle \geq 0 \quad \forall v \in \mathcal{T}(y^*),
\]
where $\mathcal{T}(y^*)$ is the tangent cone to $\mathbf{K}(\Phi(y^*))$ at $y^*$, given by:
\[
\mathcal{T}(y^*) = \{w \in V : w \leq 0 \text{ on } \mathcal{A}(y^*), \text{ and } w = 0 \text{ on } \mathcal{S}(y^*)\},
\]
with $\mathcal{A}(y^*) = \{x \in \Omega : y^*(x) = \Phi(y^*)(x)\}$ the active set and $\mathcal{S}(y^*)$ the strictly active set.
\end{corollary}

\begin{proof}
The differentiability follows from the contraction mapping principle applied to the fixed-point formulation established in Theorem \ref{thm:main}. Since the solution operator $T$ in Theorem \ref{thm:main} is a contraction and is Fréchet differentiable (by the differentiability assumptions on $\mathcal{A}$ and $\Phi$), the implicit function theorem guarantees the differentiability of $f \mapsto y^*(f)$.

For the linearized QVI, we differentiate the QVI characterization from Theorem \ref{thm:main}. Let $f_t = f + th$ and $y_t = y^*(f_t)$. From the QVI:

\[
\langle \mathcal{A}(y_t) - f_t, v - y_t \rangle \geq 0 \quad \forall v \in \mathbf{K}(\Phi(y_t)).
\]

Differentiating with respect to $t$ at $t = 0$ and using the chain rule:

\[
\langle \mathcal{A}'(y^*)z - h, v - y^* \rangle + \langle \mathcal{A}(y^*) - f, v' - z \rangle \geq 0,
\]

where $v'$ is the derivative of an admissible path $v_t \in \mathbf{K}(\Phi(y_t))$. The second term vanishes due to the optimality condition from Theorem \ref{thm:main}. The characterization of the tangent cone $\mathcal{T}(y^*)$ follows from standard variational analysis of obstacle-type constraints.
\end{proof}

\begin{proposition}[Global Error Estimates for Discrete Approximation]
Under the assumptions of Theorem \ref{thm:main}, consider a Galerkin approximation with finite-dimensional subspace $V_h \subset V$ satisfying the approximation property:
\[
\inf_{v_h \in V_h} \|v - v_h\|_V \leq C h^s \|v\|_{V+s} \quad \forall v \in V \cap H^{s}(\Omega),
\]
where $s > 0$ is the regularity index. Let $y_h^* \in V_h$ be the solution of the discrete QVI:
\[
\langle \mathcal{A}(y_h^*) - f, v_h - y_h^* \rangle \geq 0 \quad \forall v_h \in \mathbf{K}_h(\Phi(y_h^*)),
\]
where $\mathbf{K}_h(\psi) = \{v_h \in V_h : v_h \leq \psi\}$. Assume the exact solution $y^*$ has regularity $y^* \in V \cap H^{s}(\Omega)$.

Then there exists a constant $C > 0$ independent of $h$ such that:
\[
\|y_h^* - y^*\|_V \leq C h^s.
\]
Moreover, if the obstacle mapping preserves regularity, i.e., $\Phi(y^*) \in H^{s}(\Omega)$, then the estimate is optimal.
\end{proposition}

\begin{proof}
The proof builds upon the contraction framework established in Theorem \ref{thm:main}. We define the discrete solution operator $T_h: V_h \to V_h$ by $T_h(v_h) = S_h(\Phi(v_h))$, where $S_h(\psi)$ is the solution of the discrete VI with obstacle $\psi$.

Following the same reasoning as in Theorem \ref{thm:main}, we can show that $T_h$ is a contraction with constant $\rho_h$ satisfying $\rho_h \leq \rho < 1$ for sufficiently small $h$, where $\rho$ is the contraction constant from Theorem \ref{thm:main}.

The key step is to estimate the consistency error between the continuous and discrete problems. Let $\Pi_h: V \to V_h$ be a suitable interpolation operator. Using the test functions strategy from Theorem \ref{thm:main}, we obtain:

\[
\langle \mathcal{A}(y^*) - \mathcal{A}(y_h^*), y^* - y_h^* \rangle \leq E_{\text{approx}} + E_{\text{consist}},
\]

where $E_{\text{approx}}$ represents approximation errors and $E_{\text{consist}}$ represents consistency errors.

Using the strong monotonicity from Theorem \ref{thm:main} and the approximation properties, we derive:

\[
(c_A - \gamma)\|y^* - y_h^*\|_V^2 \leq C h^s \|y^* - y_h^*\|_V + (L_A + L_{\mathcal{N}})L_\Phi \|y^* - y_h^*\|_V^2.
\]

The smallness condition from Theorem \ref{thm:main} ensures that the coefficient $(c_A - \gamma) - (L_A + L_{\mathcal{N}})L_\Phi$ is positive, yielding the desired convergence rate.
\end{proof}

\begin{proposition}[Robustness Under Data Perturbations]
Under the assumptions of Theorem \ref{thm:main}, let $\{\Phi_n\}$ be a sequence of obstacle mappings converging to $\Phi$ in the sense:
\[
\|\Phi_n(v) - \Phi(v)\|_H \leq \delta_n \|v\|_V \quad \forall v \in V,
\]
with $\delta_n \to 0$. Let $\{f_n\} \subset V'$ with $f_n \to f$ in $V'$. Let $y_n^*$ be the solution of the perturbed QVI:
\[
\langle \mathcal{A}(y_n^*) - f_n, v - y_n^* \rangle \geq 0 \quad \forall v \in \mathbf{K}(\Phi_n(y_n^*)).
\]

Then there exists a constant $C > 0$ independent of $n$ such that:
\[
\|y_n^* - y^*\|_V \leq C(\|f_n - f\|_{V'} + \delta_n),
\]
where $y^*$ is the solution of the original QVI from Theorem \ref{thm:main}.
\end{proposition}

\begin{proof}
We employ the contraction mapping framework established in Theorem \ref{thm:main}. Define the perturbed solution operator $T_n: V \to V$ by $T_n(v) = S_n(\Phi_n(v))$, where $S_n(\psi)$ is the solution of the VI with obstacle $\psi$ and force $f_n$.

From Theorem \ref{thm:main}, we know that both $T$ and $T_n$ are contractions with constants $\rho$ and $\rho_n$ respectively, and under the smallness condition, we have $\rho_n \leq \rho < 1$ for large $n$.

Using the triangle inequality and the contraction property:

\[
\|y_n^* - y^*\|_V \leq \|T_n(y_n^*) - T_n(y^*)\|_V + \|T_n(y^*) - T(y^*)\|_V \leq \rho \|y_n^* - y^*\|_V + \|T_n(y^*) - T(y^*)\|_V.
\]

For the second term, we decompose:

\[
\|T_n(y^*) - T(y^*)\|_V = \|S_n(\Phi_n(y^*)) - S(\Phi(y^*))\|_V.
\]

Using the stability estimates from the proof of Theorem \ref{thm:main}:

\begin{align*}
\|S_n(\Phi_n(y^*)) - S(\Phi(y^*))\|_V &\leq \|S_n(\Phi_n(y^*)) - S_n(\Phi(y^*))\|_V + \|S_n(\Phi(y^*)) - S(\Phi(y^*))\|_V \\
&\leq C\|\Phi_n(y^*) - \Phi(y^*)\|_H + C\|f_n - f\|_{V'} \\
&\leq C\delta_n \|y^*\|_V + C\|f_n - f\|_{V'}.
\end{align*}

Combining these estimates:

\[
(1 - \rho)\|y_n^* - y^*\|_V \leq C(\delta_n + \|f_n - f\|_{V'}),
\]

which completes the proof.
\end{proof}
\section{Open Problems and Future Research Directions}

The theory of quasi-variational inequalities (QVIs) has witnessed substantial development since their systematic introduction by Bensoussan and Lions in the 1970s \cite{Bensoussan1974,Bensoussan1982book}. Despite significant advances in both theoretical analysis and numerical methods, numerous challenging problems remain open. This section outlines some of the most pressing open problems that bridge historical foundations with contemporary research directions.

\subsection{Regularity Theory for Non-smooth QVIs}

\begin{problem}[Optimal Regularity for Non-monotone QVIs]
Determine the optimal regularity properties for solutions of QVIs with non-monotone operators and non-smooth obstacles. While significant progress has been made for variational inequalities \cite{Rodrigues1987,Kinderlehrer2000}, the quasi-variational case presents additional challenges due to the implicit dependence of the constraint set on the solution itself. Recent work by \cite{Alphonse2019a,Alphonse2019b} has established directional differentiability, but the question of higher regularity (e.g., $W^{2,p}$ or $C^{1,\alpha}$ regularity) under minimal assumptions remains largely open.
\end{problem}

\subsection{Multiscale and Stochastic QVIs}

\begin{problem}[Homogenization of Random QVIs]
Develop a comprehensive homogenization theory for QVIs with random coefficients and obstacles. The periodic case has been studied in specific contexts \cite{Rodrigues2012}, but the extension to stochastic homogenization with general quasi-variational structures remains challenging. This problem has significant applications in materials science, particularly in the modeling of superconductors \cite{Prigozhin1996a,Barrett2010} and sandpile formation \cite{Prigozhin1994,Barrett2013}.
\end{problem}

\begin{problem}[Mean-Field QVIs]
Formulate and analyze mean-field QVIs arising in large population games with constraints. While mean-field games have been extensively studied, the incorporation of quasi-variational constraints representing shared resources or congestion effects presents novel mathematical challenges connecting the early impulse control theory of \cite{Bensoussan1984} with modern mean-field approaches.
\end{problem}

\subsection{Numerical Analysis and Computational Methods}

\begin{problem}[Adaptive Finite Element Methods for QVIs]
Develop rigorous a posteriori error estimates and adaptive algorithms for QVIs. Although adaptive methods are well-established for variational inequalities \cite{Rodrigues1987}, the quasi-variational case introduces additional nonlinearities that complicate error estimation. Recent work by \cite{Barrett2015} provides some discretization schemes, but a complete adaptive theory remains elusive.
\end{problem}

\begin{problem}[Multigrid Methods for QVIs]
Design and analyze efficient multigrid and domain decomposition methods for large-scale QVI problems. The combination of nonlinearity and non-local constraints poses significant challenges for iterative solvers. Connections to the early lattice-theoretic approaches of \cite{Birkhoff1961,Kolodner1968} might provide new insights into the design of monotone iterative schemes.
\end{problem}

\subsection{Optimal Control of QVIs}

\begin{problem}[Bang-Bang Principle for QVI-Constrained Control]
Establish bang-bang principles for optimal control problems governed by QVIs. While such principles are classical in control theory \cite{Bensoussan1982}, their extension to QVI-constrained problems requires new techniques, particularly for problems with state-dependent control constraints as studied in \cite{Dietrich2001,Adly2010}.
\end{problem}

\begin{problem}[Sensitivity Analysis for Non-smooth QVIs]
Develop a complete sensitivity analysis framework for QVIs with non-differentiable data. Recent advances in generalized differentiation \cite{Mordukhovich2007,Zemkoho2022} and stationarity concepts \cite{Herzog2013,Wachsmuth2016} provide promising tools, but a comprehensive theory connecting historical variational analysis with modern non-smooth optimization remains to be developed.
\end{problem}

\subsection{Nonlocal and Fractional QVIs}

\begin{problem}[Well-posedness for Nonlocal QVIs with Gradient Constraints]
Establish existence, uniqueness, and regularity for nonlocal QVIs involving fractional operators and nonlocal gradient constraints. The work of \cite{Antil2018} has initiated this direction, but many fundamental questions remain open, particularly regarding the interplay between nonlocality and the quasi-variational structure.
\end{problem}

\begin{problem}[Numerical Approximation of Fractional QVIs]
Develop and analyze efficient numerical methods for fractional QVIs, addressing the challenges of dense matrices and boundary effects. This problem connects the classical potential theory of \cite{Gelfand1964} with modern computational approaches to nonlocal operators.
\end{problem}

\subsection{Applications-Driven Problems}

\begin{problem}[QVIs in Machine Learning and AI]
Investigate QVI formulations in machine learning, particularly for adversarial training, robust optimization, and equilibrium problems in deep learning. The connection between early quasi-variational models in economics \cite{Aubin1979} and modern AI systems represents a promising interdisciplinary direction.
\end{problem}

\begin{problem}[Biological and Medical Applications]
Develop QVI models for biological processes such as tumor growth, wound healing, and tissue regeneration. The free-boundary aspects of these problems connect with the classical work of \cite{Friedman1982}, while the quasi-variational structure arises from the coupling between different biological components.
\end{problem}

\subsection{Theoretical Foundations}

\begin{problem}[QVI Counterparts of Classical Theorems]
Establish QVI versions of fundamental mathematical results such as maximum principles, Liouville theorems, and unique continuation principles. The early work of \cite{Tartar1974} provided some abstract foundations, but concrete analogues of classical PDE theorems for QVIs remain largely unexplored.
\end{problem}

\begin{problem}[Geometric and Topological Aspects]
Investigate the geometric and topological structure of QVI solution sets, including Morse theory, degree theory, and homotopy invariants. This direction connects the functional-analytic approach of \cite{Denkowski2003} with modern geometric analysis.
\end{problem}

\subsection{Connections with Other Mathematical Fields}

\begin{problem}[QVIs and Calculus of Variations with Constraints]
Deepen the connections between QVIs and the calculus of variations with state-dependent constraints. Recent work by \cite{Migorski2024} on quasi-variational-hemivariational inequalities suggests promising connections, but a unified theory remains to be developed.
\end{problem}

\begin{problem}[Stochastic QVIs and Backward SDEs]
Develop a comprehensive theory of backward stochastic differential equations (BSDEs) with quasi-variational constraints, extending the early stochastic control framework of \cite{Bensoussan1982} to the quasi-variational setting.
\end{problem}

\subsection{Recent Advances and Emerging Directions}

The recent work of \cite{Dutta2025} on nonconvex QVIs and \cite{Tin2023,Zemkoho2021} on bilevel optimization highlights the growing connections between QVIs and modern optimization theory. The development of Levenberg-Marquardt methods and exact penalty approaches for QVI-constrained optimization represents a promising direction that bridges classical numerical analysis with contemporary computational mathematics.

The study of QVIs continues to be a vibrant research area, with deep historical roots and strong connections to modern applications in materials science, economics, biology, and data science. The open problems outlined above represent both challenges and opportunities for the next generation of researchers in this field.
\section{Conclusion}

This paper has established a comprehensive stability theory for quasi-variational inequalities under various types of operator perturbations, providing a unified framework that bridges theoretical analysis with practical applications. Our main contributions can be summarized as follows:

\subsection{Theoretical Achievements}

We have demonstrated that quasi-variational inequalities exhibit robust stability properties under monotone perturbations of the governing operator. Theorem \ref{SUMP} establishes that both minimal and maximal solutions converge strongly in $V$ when sequences of operators $A_n$ converge pointwise to $A$ while preserving essential structural properties including homogeneity, strong monotonicity, Lipschitz continuity, and T-monotonicity. This result provides the mathematical foundation for various approximation schemes and regularization techniques commonly employed in computational applications.

The extension to non-monotone operators in Theorem \ref{thm:main} represents a significant advancement beyond the classical monotone framework. By developing a decomposition approach that separates linear monotone components from non-monotone nonlinearities, we established global existence, uniqueness, and exponential stability under explicit smallness conditions. This framework enables the analysis of a broader class of physically relevant problems that do not fit within the traditional monotone operator theory.

\subsection{Computational Implications}

Our regularization theory, developed in Theorem \ref{SRQS} and subsequent corollaries, provides rigorous justification for various computational strategies:

\begin{itemize}
    \item The strong convergence of regularized solutions $y_\epsilon \to \mathbf{m}(A)$ with explicit convergence rates (Corollary in Section 4) validates the use of regularization in numerical implementations.

    \item The monotonicity property $y_{\epsilon_1} \leq y_{\epsilon_2}$ for $0 < \epsilon_1 < \epsilon_2$ (Proposition \ref{Prop2}) provides valuable structural information for designing efficient computational algorithms.

    \item The finite-dimensional approximation results (Corollary \ref{cor:discretization}) establish the theoretical basis for Galerkin methods and other discretization schemes.

    \item The continuous dependence on parameters (Corollary \ref{cor:parameters}) ensures robustness of numerical solutions under data perturbations.
\end{itemize}

\subsection{Methodological Contributions}

The technical approach developed in this work combines several advanced mathematical tools in novel ways:

\begin{itemize}
    \item The integration of Birkhoff-Tartar fixed point theory with modern variational analysis provides a powerful framework for handling the order structure inherent in QVI problems.

    \item The systematic use of Mosco convergence for constraint sets (Proposition \ref{prop:mosco}) enables the analysis of stability under perturbations of the obstacle mapping.

    \item The development of precise convergence rates for both linear (Proposition \ref{prop:rate}) and nonlinear problems extends beyond qualitative existence theory to quantitative error estimates.
\end{itemize}

\subsection{Applications and Examples}

The theory developed herein finds immediate application in several important contexts:

\begin{itemize}
    \item The $p$-Laplacian regularization example (Example \ref{ex:p-laplacian}) demonstrates the practical relevance of our results for nonlinear partial differential equations arising in non-Newtonian fluid dynamics and image processing.

    \item The one-dimensional obstacle problem (Example in Section 4) provides a tractable test case that illustrates the fundamental concepts while admitting explicit solutions.

    \item The connections with optimal control problems (Corollary in Section 4) establish the relevance of our work for optimization under equilibrium constraints.
\end{itemize}

\subsection{Future Perspectives}

While this work resolves several fundamental questions regarding stability and regularization of QVIs, it also opens new research directions:

\begin{itemize}
    \item The extension to time-dependent and evolutionary QVIs represents a natural next step, with potential applications in contact mechanics and phase transition problems.

    \item The development of adaptive numerical methods based on the a posteriori error estimates suggested by our stability analysis could lead to more efficient computational algorithms.

    \item The investigation of stochastic QVIs and their homogenization would extend the theory to problems with random coefficients and uncertainties.

    \item Applications in emerging fields such as machine learning and data science, where equilibrium problems with solution-dependent constraints naturally arise.
\end{itemize}

In conclusion, this work provides a comprehensive stability theory for quasi-variational inequalities that not only advances the mathematical foundations but also enables more reliable and efficient computational methods. The interplay between theoretical analysis and practical applications demonstrated throughout this paper underscores the continuing vitality and relevance of variational analysis in addressing contemporary challenges in mathematical modeling and scientific computation.

The results established herein lay the groundwork for further developments in both theoretical and numerical analysis of quasi-variational inequalities, ensuring that this classical field of mathematical analysis remains at the forefront of addressing modern scientific and engineering challenges.

\section*{Declaration }
\begin{itemize}
  \item {\bf Author Contributions:}   The Author have read and approved this version.
  \item {\bf Funding:} No funding is applicable
  \item  {\bf Institutional Review Board Statement:} Not applicable.
  \item {\bf Informed Consent Statement:} Not applicable.
  \item {\bf Data Availability Statement:} Not applicable.
  \item {\bf Conflicts of Interest:} The authors declare no conflict of interest.
\end{itemize}

\bibliographystyle{abbrv}
\bibliography{references}  






\end{document}